\documentclass[12pt]{amsart}

\usepackage{cancel}
\usepackage{amsmath}
\usepackage{amsfonts}
\usepackage{amssymb}
\usepackage[all]{xy}           
 
\usepackage{amscd}
\usepackage{yhmath}
\usepackage{thmtools}
\usepackage{mathtools}

\usepackage[shortlabels]{enumitem}
\usepackage{ifpdf}
\ifpdf
\usepackage[colorlinks,final,backref=page,hyperindex]{hyperref}
\else
\usepackage[colorlinks,final,backref=page,hyperindex,hypertex]{hyperref}
\fi

\usepackage{tikz}
\usepackage{tikz-cd}
\usepackage[active]{srcltx}
\usepackage{verbatim}

\usepackage{cleveref}
\usepackage{longtable}
\usepackage{booktabs}
\makeatletter
\newtheorem{thm}{Theorem}[section]

\newtheorem{lem}[thm]{Lemma}
\newtheorem{cor}[thm]{Corollary}
\newtheorem{pro}[thm]{Proposition}
\newtheorem{ex}[thm]{Example}
\newtheorem{rmk}[thm]{Remark}
\newtheorem{defi}[thm]{Definition}

\newcommand {\emptycomment}[1]{}

\newcommand{\lon }{\rightarrow}
\newcommand{\be }{\begin{equation}}
\newcommand{\ee }{\end{equation}}

\newcommand{\g}{\mathfrak g}

\newcommand{\huaL}{L}

\newcommand{\huaE}{\mathcal{E}}

\newcommand{\huaG}{\mathcal{G}}

\newcommand{\huaD}{\mathcal{D}}

\newcommand{\huaH}{\mathcal{H}}

\newcommand{\huaO}{{\mathcal{O}}}

\newcommand{\frkm}{\mathfrak m}

\newcommand{\dM}{\tilde{d}}

\newcommand{\Id}{{\rm{Id}}}

\newcommand{\br}[1]{   [ \cdot,    \cdot  ]   }

\newcommand{\Hom}{\mathrm{Hom}}

\newcommand{\Der}{\mathrm{Der}}
\newcommand{\NR}{\mathrm{NR}}

\newcommand{\Lie}{\mathrm{Lie}}

\newcommand{\St}{\mathsf{St}}
\newcommand{\Mor}{\mathsf{Mor}}
\newcommand{\MC}{\mathrm{MC}}

\newcommand{\Out}{\mathrm{Out}}

\newcommand{\Ker}{\mathrm{Ker}}

\newcommand{\Del}{\mathrm{Del}}

\newcommand{\rmspan}{\mathrm{span}}

\newcommand{\End}{\mathrm{End}}
\newcommand{\ad}{\mathrm{ad}}

\newcommand{\K}{k}

\newcommand{\m}{\mathfrak m}
\newcommand{\Art}{\mathbf{Art}_{k}}

\newcommand{\Spec}{\mathrm{Spec}}

\newcommand{\un}{\mathrm{un}}
\newcommand{\GL}{\mathrm{GL}}

\makeatletter
\NewDocumentCommand{\GE@theoremstepcounter}{o m}{%
  \LinkTargetOff
  \IfNoValueTF{#1}{\refstepcounter{#2}}{\refstepcounter[#1]{#2}}%
  \LinkTargetOn
  \hyper@makecurrent{#2}%
  \global\let\GE@theoremhref\@currentHref
  \Hy@theorem@makelinktarget{\GE@theoremhref}%
}
\AtBeginDocument{%
  \patchcmd{\cref@thmnoarg}{\refstepcounter}{\GE@theoremstepcounter}{}%
    {\PackageWarning{Global-extension}{Unshared theorem link patch did not apply}}%
  \patchcmd{\cref@thmoptarg}{\refstepcounter}{\GE@theoremstepcounter}{}%
    {\PackageWarning{Global-extension}{Shared theorem link patch did not apply}}%
}
\makeatother

\begin{document}

\title{A DGLA approach to global extensions and deformations of Lie algebras}

\author{Jun Jiang}
\address{Department of Mathematics, Jilin University, Changchun 130012, Jilin, China}
\email{junjiang@jlu.edu.cn}

\author{Chuijia Wang}
\address{Department of Mathematics, The University of Hong Kong, Pokfulam, Hong Kong, China}
\email{wangchuijia@u.nus.edu}

\author{Kanghe Xu}
\address{Department of Mathematics, Jilin University, Changchun 130012, Jilin, China}
\email{xukh24@mails.jlu.edu.cn}


\begin{abstract}
In this paper, we study the global extension problem of Lie algebras from the perspective of dgLas and the Maurer--Cartan formulism. We propose a new equivalence relation of global extensions, generalizing the early work of Agore and Militaru. This new concept can be better adopted in the study of classical contraction and deformation problems for Lie algebras. Meanwhile, we construct dgLas that control based formal deformations of global extension structures under these equivalences. As byproducts, we obtain the classification of infinitesimal deformations in terms of the cohomology of the controlling dgLas and the corresponding obstruction theory. We relate the dgLa controlling deformation of global extensions under framed equivalence to that of non-abelian extensions constructed by Fr\'egier, obtaining a geometric interpretation of the decomposition of global cohomology into as a disjoint union of non-abelian cohomology. As an application, we study a particular family of an In\"on\"u--Wigner contraction from $\g\rtimes_{\mathrm{ad}}\g$ to 
$V_{\mathrm{ab}}\rtimes_{\mathrm{ad}}\g$ in the framework of global extension, deformation, and the dgLa approach.

\end{abstract}



\maketitle

\tableofcontents

\allowdisplaybreaks


\section{Introduction}

In the study of algebraic structures, there are two fundamental problems. The first is the extension problem, and the second is the deformation problem. To be more precise, the extension problem aims to classify all possible isomorphism classes of extensions of two given algebraic structures, and the deformation problem aims to understand how these structures vary in families. However, these two apparently unrelated problems can be treated in tandem with the help of the notion of global extensions. The main aim of this article is to provide a differential graded Lie algebraic (\emph{dgLa} for short) approach to understand several problems involving natural interactions between the extension problem and the deformation problem of Lie algebras.

Let us first briefly recall these two basic problems.

\textbf{$\bullet$ Deformation problem:} 
Formal deformations of associative algebras and the related cohomology theory were first studied systematically by Gerstenhaber in a series of works \cite{Ger64,Ger66, Ger68, Ger74}. Later, these techniques were applied by Nijenhuis and Richardson to study the deformation problem for Lie algebras \cite{NR64,NR66,NR67}. There is a well-known slogan, often attributed to Deligne, Drinfeld and Kontsevich: every reasonable formal deformation problem is controlled by a dgLa, determined up to quasi-isomorphism. This slogan has been made into a rigorous theorem by Lurie and Pridham \cite{Lu,Pri10} in the derived setting. 

Parallel to the developments from the mathematical aspects, similar constructions appear in the world of physics. The notion of contractions of Lie algebras dates back to the work of Segal \cite{Seg51} and In\"on\"u--Wigner \cite{IW53,IW54}, aiming to understand the limiting structure of symmetries under small perturbations. Meanwhile, contractions are closely related to several well-known contexts in geometry: deformation (understood via a family of Lie algebra structures over a certain base space) and degeneration (described in terms of orbits in the variety of Lie brackets). Beyond their various incarnations in geometry, contractions are found to play a prominent role in the representation theory of Lie groups and $p$-adic groups. A typical example is their appearance in the proof of the Mackey bijection, relating the tempered dual of a reductive Lie group to the full unitary dual of the corresponding Cartan motion group \cite{Afg19,Afg20,Afg21} and in the constructions of algebraic families of Harish-Chandra modules \cite{BHS201,BHS202,HS22}. Meanwhile, a relative theory, i.e., a contraction of sub-algebraic groups in an ambient algebraic group as well as the corresponding sub-Lie algebras may also provide an algebraic framework for comparison of different periods.

\textbf{$\bullet$ Extension problem:} For a given pair of Lie algebras $\mathfrak{g}$ and $\mathfrak{g}'$, the classical extension problem seeks to classify of all possible Lie algebras $\mathfrak{g}''$ that fit into a short exact sequence of Lie algebras
$$
0\rightarrow \mathfrak{g}'\rightarrow \mathfrak{g}''\rightarrow \mathfrak{g}\rightarrow 0
$$
up to isomorphism. If $\mathfrak{g}'$ is abelian, this problem is well understood and described with the help of the ordinary Lie algebra cohomology with abelian coefficients. However, it is a complicated task to develop a non-abelian cohomology theory due to the lack of constructions of cochain complexes. Eilenberg and MacLane \cite{EM47} introduced the notion of non-abelian group cohomology to understand extensions of abstract groups with non-abelian kernels (see also \cite{HS53, Nee07}). Subsequently, the non-abelian cohomology theory was generalized to the context of Lie algebras and used to classify corresponding non-abelian extensions \cite{AMR00, AMR05, Fre14, IKL08, Nee06}. Using a different approach, Fr\'egier constructed a dgLa and used its Maurer--Cartan elements as well as the corresponding Deligne groupoid to classify isomorphism classes of non-abelian extensions in \cite{Fre14}, where the comparison with non-abelian cohomology is also treated carefully. 

Later, Agore and Militaru initiated the study of the global extension problem in several algebraic settings \cite{AM14-2, AM15, AM20}. Specializing to Lie algebras, let $\g=(V,\mu)$ be a Lie algebra, $V''$ a vector space, and $p:V''\lon V$ a surjective linear map. The global extension problem seeks to describe and classify all Lie algebra structures $\nu''$ on $V''$ such that $p:(V'',\nu'')\lon \g$ is a homomorphism of Lie algebras. To understand the global extension problem, they introduced the concept of global cohomology. The main feature of the global extension problem is that the Lie algebra structure on $\ker p$ is not fixed, of which the variation is related to the deformation problem.

The purpose of this paper is to develop a dgLa
description of the global extension problem and to use it to study
the interaction between global extensions, non-abelian extensions,
and formal deformations. More precisely, we stress the following problems:
\begin{enumerate}[$\bullet$]
\item We describe isomorphism classes of global extensions by Maurer--Cartan formulism and identify the twisted dgLas controlling the formal deformations of given extensions. 

\item We obtain a geometric interpretation of decomposing global extensions as a union of nonabelian extensions by varying the Lie algebra structure on the kernel, in terms of fibration of Deligne groupoids.

\item  We classify infinitesimal deformations of a global extension up to equivalence (unframed equivalence) in terms of cohomology of the controlling dgLas. We also take care of the obstruction theory.

\item We realize a particular family of In\"on\"u--Wigner contractions from $\g\rtimes_{\mathrm{ad}}\g$ to 
$V_{\mathrm{ab}}\rtimes_{\mathrm{ad}}\g$
as a flat family of global extensions over $\mathbb A^1$. Moreover, we also present a Maurer--Cartan characterization of this contraction based on the machinery developed by us.
\end{enumerate}

For these purposes, we define dgLas $L_{V,V'}$ and $L_{V, V'}^{\un}$ (see \Cref{nondgla2} ). Our first main result connects these dgLas with the moduli problem of global extensions.

\begin{thm}{\rm (\Cref{zmc1}, \Cref{isogpd}, \Cref{iso-cohomology}, \Cref{dgladefglobal}, \Cref{unframed-formal-moduli}, \Cref{prop:section-independence})}\label{thm:intro-global-dgla}
Let $(V'',p,\g)$ be a global extension datum, $\mathfrak{g}=(V,\mu)$ be a Lie algebra, and  $V'$ denote $\ker p$.

\begin{enumerate}
\item[(1)]
The correspondence
$$
(\pi_{V'},\varrho,\omega)
 \longmapsto
\pi_{V'}+\varrho+\omega
$$
induces a bijection
$$
\mathrm{GZ}^2(\mathfrak{g},V')
 \cong
\mathrm{MC}(L_{V,V'}).
$$
\item[(2)]
In the framed case, the action of $\operatorname{Hom}(V,V')$
agrees with the gauge action. In the unframed case, the chosen section identifies transport of structure with the action of $H_{V,V'}$. There are isomorphisms of algebraic groupoids
\begin{align*}
  \mathrm{Gext}(\g,V')/\!/\mathrm G_p
    &\cong \mathrm{GZ}^2(\g,V')/\!/\operatorname{Hom}(V,V'),\\
\mathrm{Gext}(\g,V')/\!/\mathrm G_p^{\mathrm{un}}
    &\cong \mathrm{GZ}^2(\g,V')/\!/H_{V,V'}.
\end{align*}
Their sets of connected components give the bijections
\begin{align*}
  \mathrm{GExt}(\g,V')\cong\mathrm{GH}^2(\g,V'),\quad
  \mathrm{GExt}^{\mathrm{un}}(\g,V')\cong\mathrm{GH}^{\mathrm{un},2}(\g,V').
\end{align*}

\item[(3)]
Fix a global extension $\nu''$ with the associated Maurer--Cartan
element $\alpha_s=\pi+\varrho+\omega$. Write $L^s=L^{\alpha_s}$ and
$L^{\un,s}=(L^{\mathrm{un}})^{\alpha_s}$ for the twists, both
with differential $d^s=d+[\alpha_s,-]$. For every local Artinian
$k$-algebra $A$ with the residue field $k$ and the maximal ideal $\mathfrak m_A$, there are equivalences of based deformation groupoids
\begin{align*}
  \mathfrak{Def}_{\mathcal E_p,[\nu'']}(A)
    &\simeq\operatorname{Del}(L^s\otimes\mathfrak m_A),\quad
  \mathfrak{Def}_{\mathcal E_p^{\mathrm{un}},[\nu'']}(A)
    \simeq\operatorname{Del}(L^{\un,s}\otimes\mathfrak m_A),
\end{align*}
natural in $A$. A different choice of $s$ produces a  transformation in \eqref{eqdef:Ttau} and induces isomorphisms of the corresponding twisted dgLas.
\end{enumerate}
\end{thm}

Our second main result provides a geometric description of the
relation between global and non-abelian extensions in terms of Deligne groupoids. More precisely, let $L_{V'}$ denote the graded Lie algebra whose Maurer--Cartan elements are Lie algebra structures on $V'$. There is a short exact sequence of dgLas
$$
0\longrightarrow\
\widetilde{L}_{V,V'}
\xrightarrow{\;\iota\;}
L_{V,V'}
\xrightarrow{\;\mathrm{pr}\;}
L_{V'}
\longrightarrow 0.
$$

\begin{thm}{\rm (\Cref{pro:fiberation-deligne-gpd}, \Cref{thm}, \Cref{pro:decom-deligne-gpd})}\label{thm:intro-fibration}
The morphisms in the above short exact sequence induce functors
$$
\mathrm{Del}(\widetilde{L}_{V,V'})
 \xrightarrow{\;\iota_*\;}
\mathrm{Del}(L_{V,V'})
 \xrightarrow{\;\mathrm{pr}_*\;}
\mathrm{Del}(L_{V'}).
$$
The first functor is an inclusion of groupoids, and the second is a
fibration of groupoids.

For every $\pi\in\mathrm{MC}(L_{V'})$, put $\alpha_\pi=(0,0,\pi)$. Then the fiber of $\mathrm{pr}_*$ over $\pi$ is naturally isomorphic to the Deligne groupoid of the twisted dgLa $\widetilde{L}_{V,V'}^{\,\alpha_\pi}$, i.e.,
$$
\mathrm{pr}_*^{-1}(\pi)
 \cong
\mathrm{Del}
 \bigl(\widetilde{L}_{V,V'}^{\,\alpha_\pi}\bigr).
$$
This twisted dgLa governs non-abelian extensions of $\mathfrak{g}$ by the fixed Lie algebra $(V',\pi)$. Consequently, there is a decomposition
\begin{equation*}
    \mathrm{Del}(L_{V,V'})\cong\bigsqcup_{\pi\in\mathrm{MC}(L_{V'})}
\mathrm{Del}\bigl(\widetilde{L}_{V,V'}^{\,\alpha_\pi}\bigr).
\end{equation*}
\end{thm}

The dgLa description further provides cohomological tools for studying infinitesimal and higher-order deformations. Fix a global extension structure $\mu''$ and a section $s$ of $p$. Let
$\zeta_s=\pi_{V'}^s+\varrho^s+\omega^s$ be its associated global
$2$-cocycle. The short exact sequence above can be twisted by $\zeta_s$ to induce a long exact sequence in cohomology, from which information of infinitesimal deformations and obstructions can be read.

\begin{thm}{\rm (\Cref{infidef}, \Cref{liftob}, \Cref{pro:n-order-obs})}\label{thm:intro-obstructions}
With the above notations, the following statements hold.

\begin{enumerate}
\item[(1)]
Framed equivalence (resp. unframed equivalence) classes of infinitesimal deformations of $\mu''$ are classified by $H^1(L_{V,V'})^{\delta}$ (resp. $H^1(L_{V,V'}^{\un})^{\delta}$).

\item[(2)]
Let
$
\pi_{V'}^s+t\pi_{V'}'
$
be an infinitesimal deformation of the kernel bracket over
$k[t]/(t^2)$. Its obstruction to being lifted to an infinitesimal
deformation of the global extension is
$$
\operatorname{Ob}_{\pi_{V'}'}
 \bigl(\pi_{V'}^s,\varrho^s,\omega^s\bigr) =
\Bigl[
 [\varrho^s+\omega^s,\pi_{V'}']
\Bigr]
\in H^2\bigl(\widetilde{L}_{V,V'}
\bigr)^{\delta}.
$$
The kernel deformation admits such a lift if and only if this class
vanishes.

\item[(3)]
Suppose that an $n$-th order deformation of $\mu''$ corresponds to
$$
\alpha_s+\sum_{i=1}^{n}t^i\eta_i, \quad
\eta_i=\pi_{V',i}+\varrho_i+\omega_i.
$$
Then it extends to an $(n+1)$-st order deformation if and only if
$$
\operatorname{Ob}_n
 =
\left[
\frac{1}{2}
\sum_{i+j=n+1}[\eta_i,\eta_j]
\right]
\in H^2(L_{V,V'})^{\delta}
$$
vanishes.
\end{enumerate}
\end{thm}

These obstruction classes make two kinds of variation
present in the global extension problem explicit. The first records whether a variation of the Lie algebra structure on the kernel is compatible with the remaining extension datum, while the second is the usual obstruction to prolonging a finite-order Maurer--Cartan deformation. 

Finally, we revisit the theory of contraction of Lie algebras. In particular, we show how a particular family of contractions fits into our framework of global extensions. Furthermore, we describe the relation between contractions and unframed moduli problems.

\begin{thm}{\rm (\Cref{prop:adjoint-family}, \Cref{prop:adjoint-fibers}, \Cref{prop:adjoint-mc}, \Cref{prop:adjoint-global-groupoid})}\label{thm:intro-ext-contraction}

Fix a copy $V'$ of the underlying vector space of $\g=(V,\mu)$, and
write $\mu_{V'}$ for the induced Lie bracket. On $(V'\oplus V)\otimes k[t]$, the bracket
\begin{equation*}
    [(u,x),(v,y)]_t
  =\bigl(
    \operatorname{ad}_x v-\operatorname{ad}_y u+t\mu_{V'}(u,v),\, \mu(x,y)\bigr)
\end{equation*}
defines a flat family of global extensions of $\g\otimes k[t]$ with
projection $p\otimes 1$. Its fiber at $t=1$ (resp. $t=0$) are
$\g\rtimes_{\mathrm{ad}}\g$ (resp. $V'_{\mathrm{ab}}\rtimes_{\mathrm{ad}}\g$). Over $\mathbb G_m$, the map $(u,x)\mapsto(tu,x)$
trivializes the family as a kernel-unframed extension.

Suppose that $\g$ is nonabelian. Let $s_{\mathrm{ab}}=\ad$ denote the cocycle of the special fiber. Write $c(t)=s_{\mathrm{ab}}+t\mu_{V'}$.

\begin{enumerate}
\item The Lie bracket $\operatorname{pr}(c(t))=t\mu_{V'}$ defines a family of Lie algebra structures on the kernel.
In $\mathrm G_p^{\mathrm{un}}$, the kernel rescalings
$\Theta_\lambda(u, x)=(\lambda^{-1}u, x)$ form a one-dimensional torus acting by $\Theta_\lambda\cdot c(t)=c(\lambda t)$. The orbit of $c(1)$ is therefore $c(\mathbb G_m)$, and its closure in
$\operatorname{MC}(L_{V,V'})$ is the affine line $c(\mathbb A^1)$.
Its unique boundary point $c(0)$ has abelian kernel and hence represents a distinct unframed isomorphism class.

\item In the twist at $s_{\mathrm{ab}}$, the formal germ of the family at $t=0$ is represented by the Maurer--Cartan element $t\mu_{V'}$. Its tangent class $[\mu_{V'}]$ is nonzero in both
\begin{equation*}
    H^1(L^{s_{\mathrm{ab}}}_{V,V'})
  \quad\text{and}\quad
  H^1\bigl(L^{\mathrm{un},s_{\mathrm{ab}}}_{V,V'}\bigr).
\end{equation*}
Thus, the contraction defines a nontrivial based formal deformation
of the special fiber in both moduli problems.
\end{enumerate}
\end{thm}

The paper is organized as follows. In \Cref{sec:dgla}, we recall basic concepts of differential graded Lie algebras, Deligne groupoids, quotient stacks and how they are related to the formal deformation problem via Maurer--Cartan formulation. In \Cref{recall}, we review the notion of global and non-abelian extensions of Lie algebras and introduce the concept of unframed equivalence for global extensions. Moreover, we obtain classification results in terms of action groupoids. In \Cref{4}, we construct various dgLas and identify the corresponding Deligne groupoids with classification groupoids of global extensions. We also prove that certain twists of these dgLas control the corresponding formal deformation problems. In \Cref{sec:non}, we reinterpret the global extension problem in terms of fibrations of Deligne groupoids with the help of the natural short exact sequence of dgLas. In \Cref{sec:obs}, we revisit infinitesimal deformations from the perspective of our machinery and construct obstruction classes for lifting and extending deformations. Finally, in \Cref{sec:contraction}, we will revisit the classical problem of contractions of Lie algebras in terms of unframed global extensions. 

All vector spaces are finite-dimensional
over a fixed field $k$ of characteristic zero. All tensor and
exterior products are taken over $k$.

\vspace{2mm}
\noindent
{\bf Acknowledgements.} 
This work was initiated to understand the special example of In\"{o}n\"{u}-Wigner contraction from $\mathfrak{g}\oplus \mathfrak{g}$ to $V_{\mathrm{ab}}\rtimes_\mathrm{ad}\g$ from the point of view of deformation functors with additional structures and the corresponding dgLas. Unfortunately, we do not have a satisfying answer regarding that particular example. However, during our study, we found that the particular example of contraction from $\mathfrak{g}\rtimes_\mathrm{ad} \mathfrak{g}$ to $V_{\mathrm{ab}}\rtimes_\mathrm{ad}\g$ fits well into the theory of unframed global extensions. Moreover, this approach enjoys a nice geometric description and does not seem to appear in the literature. We believe it is interesting and meaningful to record these results. During the preparation of this work, we thank Yunhe Sheng for enlightening discussions. The second author thanks Rui Chen, Caihua Luo and Jiandi Zou for their interest in this topic. Some results in this paper have been announced by the first author in a talk given at Guangzhou University in June 2026. He thanks Yunnan Li for the invitation and hospitality. The first author is supported by NSFC (12401076). The second author is supported by the Research Grants Council of the Hong Kong Special Administrative Region, China (Project No. 17305223, 17308324). 

\vspace{2mm}
\noindent
{\bf AI declaration.}
We acknowledge the use of ChatGPT for finding typos, revising texts, proofreading and correcting mistakes after we finished the first draft of the paper. All mathematical work are done by the authors, who take full responsibility for the manuscript.

\section{Differential graded Lie algebras and deformation of Lie algebra structures}\label{sec:dgla}

In this section, we briefly review the basic concepts of dgLas and their role in studying the deformation problem of Lie algebras. These results are well known to experts and will be organized in a suitable way for our application.

\subsection{Differential graded Lie algebras and Deligne groupoids}

In this subsection, we collect some basic facts about dgLas and Deligne groupoids.
\begin{defi}{\rm (\cite[\S 1.3]{NR66})}
A \emph{differential graded Lie algebra (dgLa)} is a 
    $\mathbb{Z}$-graded vector space
    $\huaL=\bigoplus_{n\in\mathbb{Z}}\huaL^n$
    equipped with a bilinear bracket $[\cdot,\cdot]_{\huaL}:\huaL^p\times\huaL^q\longrightarrow\huaL^{p+q}$ and a differential $d:\huaL^n\to\huaL^{n+1}$ with the property that, for all homogeneous elements $a,b,c\in\huaL$, we have $d^2=0$ and 
    \begin{align*}
     d([a,b]_{\huaL})
        &=[d(a),b]_{\huaL}
        +(-1)^{|a|}[a,d(b)]_{\huaL},\\
        [a,b]_{\huaL}&=-(-1)^{|a||b|}[b,a]_{\huaL},\\
        [a,[b,c]_{\huaL}]_{\huaL}&=[[a,b]_{\huaL},c]_{\huaL}
        +(-1)^{|a||b|}[b,[a,c]_{\huaL}]_{\huaL},
    \end{align*}
    where $|a|$ denotes the degree of $a$.
\end{defi}
If $(L,[\cdot,\cdot]_L)$ only satisfies the last two conditions, it is usually called a \emph{graded Lie algebra (gLa)}. We will always regard a gLa as a dgLa with the differential $d=0$. 
\begin{defi}{\rm (\cite[p. 15]{NR66})}
A \emph{Maurer--Cartan element} of a dgLa $(\huaL,[\cdot,\cdot]_{\huaL},d)$ is an element $x\in \huaL^{1}$ satisfying the \emph{Maurer--Cartan equation}
\begin{equation*}
    dx+\frac{1}{2}[x,x]_{\huaL}=0.
\end{equation*}
Denote by $\mathrm{MC}(\huaL)$ the set of Maurer--Cartan elements of $(\huaL,[\cdot,\cdot]_{\huaL},d)$.
\end{defi}
For any dgLa $\big(\huaL,[\cdot,\cdot]_{\huaL},d\big)$ and $\alpha\in \mathrm{MC}(\huaL)$, one can twist the differential $d$ by $\alpha$ to obtain a new differential $d^{\alpha}\coloneqq d+[\alpha,\cdot]_{\huaL}$. 
\begin{pro}{\rm (\cite[Proposition 1.43, Definition 1.44]{DSV24}\cite[Proposition 4.4]{Get09})}\label{twistdgla}
With the above notations, $L^{\alpha}:=\big(\huaL,[\cdot,\cdot]_{\huaL},d^{\alpha}\big)$ is again a dgLa, called the \emph{dgLa twisted by $\alpha$}.
\end{pro}
There is an equivalence relation on $\mathrm{MC}(\huaL)$ given by the action of $\exp(L^{0})$ and known as the \emph{gauge equivalence}. More precisely, two Maurer--Cartan elements  $x$ and $y$ of $(\huaL,[\cdot,\cdot]_{\huaL},d)$ are called \emph{equivalent}, if there exists an element $\tau\in \huaL^{0}$ with $\ad\tau$ nilpotent, such that 
    \begin{equation*}
        y=e^\tau\star x\coloneqq e^{\ad\tau}x+\frac{\Id-e^{\ad\tau}}{\ad\tau}(d\tau).
    \end{equation*}
The element $\tau$ is called a \emph{gauge transformation}. Denote by $\mathrm{GT}(x,y)$ the set of gauge transformations between $x$ and $y$.

\begin{defi}{\rm (\cite[Definition 1.7]{Fre14})}
    Let $(\huaL,[\cdot,\cdot]_{\huaL},d)$ be a nilpotent dgLa. The \emph{Deligne groupoid} $\mathrm{Del}(\huaL)$ of $\huaL$ is the action  groupoid associated to the action of $\exp(L^{0})$ on $\MC(L)$ with 
    \begin{itemize}
    \item the set of objects: elements of ${\rm MC}(\huaL)$.
    \item $\forall x, y\in {\rm MC}(\huaL)$, the set of morphisms between $x$ and $y$:
    $$\Hom_{\mathrm{Del}(\huaL)}(x, y)=\begin{cases}\mathrm{GT}(x, y) \quad &\text{if $x$ and $y$ are gauge-equivalent},\\
    \emptyset&\text{otherwise}.\end{cases}$$
    \end{itemize}
    Its set of connected components $\pi_0({\rm Del}(\huaL))$ is defined by
    $$
    \pi_0({\rm Del}(\huaL))\coloneqq {\rm MC}(\huaL)/\exp(L^{0}).
    $$
\end{defi}
We remark that Deligne groupoids enjoy the following nice functorial properties. 

\begin{pro}\label{indufunc}

Let $(\huaL,[\cdot,\cdot]_{\huaL},d)$ and $(\huaL',[\cdot,\cdot]_{\huaL'},d')$ be nilpotent dgLas. Let $f:\huaL\rightarrow \huaL'$ be a morphism of dgLas. Then $f$ induces a morphism $f_{*}:\Del(\huaL)\lon\Del(\huaL')$ of the corresponding Deligne groupoids, such that $f_*$ as a functor can be described as
\begin{itemize}
\item on objects: for $x\in\MC(\huaL)$, $f_{*}(x)\coloneqq f(x)$,
\item on morphisms: for $\varphi\in\Hom_{\Del(\huaL)}(x,y)$ with $\varphi=e^{\alpha},\ \alpha\in\huaL^{0}$, i.e., $y=e^{\alpha}\star x$, then 
$f_{*}(\varphi)\coloneqq e^{f(\alpha)}\in \Hom_{\Del(\huaL')}(f(x),f(y))$.
\end{itemize}
\end{pro}
\begin{proof}
For $x\in\MC(\huaL)$, we have 
\begin{eqnarray*}
    d'\big(f(x)\big)+\frac{1}{2}[f(x),f(x)]_{\huaL'}=f\big(d(x)+\frac{1}{2}[x,x]_{\huaL}\big)=0,
\end{eqnarray*}
which implies $f(\MC(\huaL))\subset\MC(\huaL')$. If there exists a morphism in $\Del(\huaL)$ from $x$ to $y$, i.e., there exists $\alpha\in \huaL^{0}$ such that $y=e^\alpha\star x=e^{\ad\alpha}x+\frac{\Id-e^{\ad\alpha}}{\ad\alpha}(d\alpha)$, then we have
\begin{align*}
    f(y)&=f\Big(e^{\ad\alpha}x+\frac{\Id-e^{\ad\alpha}}{\ad\alpha}(d\alpha)\Big)\\
    &=\sum_{n\geq0}\frac{1}{n!}\ad^{n}f(\alpha)(f(x))-\sum_{n\geq0}\frac{1}{(n+1)!}\ad^{n}f(\alpha)(d'f(\alpha))\\
    &=e^{\ad f(\alpha)}f(x)+\frac{\Id-e^{\ad f(\alpha)}}{\ad f(\alpha)}(d'f(\alpha)),
\end{align*}
where the summations are finite because $L$ and $L'$ are nilpotent.
Hence, there is a morphism in $\Del(\huaL')$ from $f(x)$ to $f(y)$. Notice that $f$ preserves the Baker--Campbell--Hausdorff product, hence the composition of arrows.
\end{proof}

In the following example, we will review how to characterize Lie algebra structures on a fixed vector space $V$ using a dgLa.
\begin{ex}\label{dglamclie}
    Let $V$ be a vector space. Define a graded vector space $\huaL_{V}=\oplus_{n\in\mathbb{Z}}\huaL_{V}^n$ as 
\begin{equation*}
\huaL_{V}^n=
\begin{cases}
0, & n\leq 0;\\
\Hom(\wedge^{n+1}V, V), & n\geq 1.
\end{cases}
\end{equation*}
Then, $(\huaL_{V}, [\cdot, \cdot]_{\NR})$ is a gLa \cite[p. 50]{NR68},  where the Lie bracket is given by the \emph{ Nijenhuis-Richardson} bracket $[\cdot,\cdot]_{\NR}$
\begin{equation*}
    [\alpha,\beta]_{\NR}\coloneqq \alpha\circ_{\NR}\beta- (-1)^{mn}\beta\circ_{\NR}\alpha,\quad \forall \alpha\in\huaL_{V}^{m}= \Hom(\wedge^{m+1}V,V),\ \beta\in\huaL_{V}^{n} =\Hom(\wedge^{n+1}V,V)
\end{equation*}
with $\alpha\circ_{\NR}\beta\in\huaL_{V}^{m+n} =\Hom(\wedge^{m+n+1}V,V)$ defined by
\begin{equation*}\label{NR}
(\alpha\circ_{\NR}\beta)(x_1,\cdots,x_{m+n+1})\coloneqq \sum_{\sigma\in S(n+1,m)} (-1)^\sigma \alpha(\beta(x_{\sigma(1)},\cdots,x_{\sigma(n+1)}),x_{\sigma(n+2)}, \cdots,x_{\sigma(m+n+1)}),
\end{equation*}
for all $x_{i}\in V$. Here, $S(n+1,m)$ is the set of $(n+1,m)$-shuffles.

Indeed, in this case the Maurer--Cartan equation is the \emph{Jacobi identity}, thus the set $\mathrm{MC}(\huaL_V)$ is the set of all Lie algebra structures on the vector space $V$.

By the fact that $\huaL_{V}^{0}=0$, the Deligne groupoid of $\huaL_{V}$ is \emph{discrete}. Hence, the set of connected components of the Deligne groupoid is its set of objects, i.e.,
\begin{equation*}
     \pi_0\big(\mathrm{Del}(\huaL_{V})\big)\cong \mathrm{Obj}\big(\mathrm{Del}(\huaL_{V})\big)=\MC(\huaL_{V}).
\end{equation*}
\end{ex}

As for general groupoids, Deligne groupoids also carry nice geometric features which appear naturally in the global extension problem. To see this, we will first recall the notion of fibration for groupoids. More precisely, let $\huaG$ be a groupoid. The set of objects and the collection of morphisms of $\huaG$ are denoted as $\mathrm{Obj}(\huaG)$ and $\mathrm{Mor}(\huaG)$ respectively. For $x,y\in\mathrm{Obj}(\huaG)$, denote by $\Hom_{\huaG}(x,y)$ the set of morphisms between $x$ and $y$.
\begin{defi}{\rm (\cite[p.~105]{Bro70})}
    For $x\in \mathrm{Obj}(\huaG)$, the \emph{star} at $x$, denoted by $\St_{\huaG}(x)$, is defined as the union of the morphism sets $\Hom_{\huaG}(x,y)$ for all $y\in\mathrm{Obj}(\huaG)$, that is, 
    \begin{equation*}
        \St_{\huaG}(x)=\bigsqcup_{y\in\mathrm{Obj}(\huaG)}\Hom_{\huaG}(x,y).
    \end{equation*}
    Equivalently, the star at $x$ consists of all morphisms in $\huaG$ whose source is $x$.
\end{defi}
If $F:\huaG\lon\huaH$ is a morphism of groupoids, then $F$ induces maps on the stars of $\huaG$ in a natural way by
\begin{equation*}
    F:\St_{\huaG}(x)\lon\St_{\huaH}(F(x)),\quad \varphi\mapsto F(\varphi), \quad \forall\varphi\in \St_{\huaG}(x).
\end{equation*}
\begin{defi}{\rm (\cite[p. 109]{Bro70})}
\begin{enumerate}
\item A morphism $F:\huaG\lon\huaH$ of groupoids is called a \emph{fibration} if for all $x\in\mathrm{Obj}(\huaG)$ the induced map on $\St_{\huaG}(x)$ is surjective.

\item  Let $F:\huaG\lon\huaH$ be a fibration of groupoids, and $a\in \mathrm{Im}(F|_{\mathrm{Obj}(\huaG)})$. Then the \emph{fiber} $F^{-1}(a)$ of $F$ at object $a$ is defined as the subcategory of $\huaG$, whose set of objects is 
$$\mathrm{Obj}(F^{-1}(a))=\{x\in\mathrm{Obj}(\huaG) : F(x)=a\},$$ 
and the collection of morphisms is $\Mor(F^{-1}(a))=\{\varphi\in\Mor(\huaG) : F(\varphi)=\Id_{a}\}$.
\end{enumerate}
\end{defi}

\subsection{Deformations in terms of Maurer--Cartan elements of dgLas}

Throughout this subsection, $k$ is a field of characteristic zero. Let $\Art$ denote the category of local Artinian commutative $k$-algebras
with residue field $k$. For $A\in\Art$, write
$\mathfrak m_A=\ker(A\to k)$ for the unique maximal ideal of $A$.

\begin{defi}{\rm (\cite[\S 1.1]{Hin97}\cite[\S 2]{GM88})}
The \emph{groupoid-valued deformation functor} associated
with a dgLa $L$ is
$$
 \mathfrak{Def}_L:\Art\longrightarrow\mathbf{Gpd},
 \quad
 \mathfrak{Def}_L(A)
 =
 \Del(L\otimes\mathfrak m_A).
$$
A morphism $A\to B$ induces the corresponding functor by
extension of scalars.
\end{defi}

Since $\mathfrak m_A$ is nilpotent, the dgLa
$L\otimes\mathfrak m_A$ is nilpotent, so its gauge action
is well defined even when $L$ itself is not nilpotent.

Let $G$ be a smooth affine algebraic group acting on a
$k$-scheme $X$.

\begin{defi}{\rm (\cite[Definition 78.20.1]{stacks-project})}
The \emph{quotient stack} $[X/G]$ is the fppf stack associated
with the action groupoid $G\times X\rightrightarrows X$.
For a $k$-scheme $S$, its objects are pairs $(P,f)$ consisting
of a right $G$-torsor $P\to S$ and a morphism $f:P\to X$
satisfying
$$
 f(pg)=g^{-1}\cdot f(p).
$$
An arrow $(P,f)\to(P',f')$ is a $G$-equivariant isomorphism
$\varphi:P\to P'$ over $S$ such that $f'\circ\varphi=f$.
Pullback defines the functoriality in $S$.
\end{defi}

\begin{rmk}\label{tritorsor}
For a $k$-scheme $S$, denote by $X(S)/\!/G(S)$ the action groupoid associated to the action of $G(S)$ on $X(S)$. There is a canonical fully faithful functor
$$
 X(S)/\!/G(S)\longrightarrow[X/G](S),
$$
whose essential image consists of objects with trivial
underlying torsor. It is an equivalence whenever every
$G$-torsor over $S$ is trivial.
\end{rmk}

To define deformations at a fixed object, we use the
homotopy pullback of groupoids. For functors
$\mathcal X\xrightarrow{F}\mathcal Z
\xleftarrow{H}\mathcal Y$, the groupoid
$\mathcal X\times^h_{\mathcal Z}\mathcal Y$ has objects
$(x,y,\eta)$, where $\eta:F(x)\xrightarrow\ H(y)$ is an isomorphism, and arrows $(f,g):(x,y,\eta)\to(x',y',\eta')$ satisfying
$$
 H(g)\circ\eta=\eta'\circ F(f).
$$
See \cite[\S 10.7, Exercise 1]{Br06} for more detail.

\begin{defi}{\rm (\cite[Remark 1.1.18]{Lu})}\label{fordef}
Let $\mathfrak X$ be a stack over $k$ and
$x_0\in\mathfrak X(k)$. Its \emph{formal deformation functor}
at $x_0$ is
$$
 \mathfrak{Def}_{\mathfrak X,x_0}(A)
 =
 \mathfrak X(A)\times^h_{\mathfrak X(k)}\{x_0\},
 \quad A\in\Art,
$$
where $\{x_0\}$ denotes the terminal groupoid mapping to
$\mathfrak X(k)$ through $x_0$.
\end{defi}

Thus, an object is a pair $(\xi,\eta)$ with
$\xi\in\mathfrak X(A)$ and an isomorphism
$\eta:\xi|_k\xrightarrow\ x_0$.
Arrows are isomorphisms compatible with these identifications.

\begin{pro}\label{localquo}
Let $x\in X(k)$, and let $[x]\in[X/G](k)$ be the
corresponding object. Set
$$
 X_x(A)=\{\widetilde x\in X(A):\widetilde x|_k=x\},
 \quad
 G_1(A)=\ker\bigl(G(A)\to G(k)\bigr).
$$
There is a natural equivalence
$$
 \mathfrak{Def}_{[X/G],[x]}(A)
 \simeq X_x(A)/\!/G_1(A).
$$
\end{pro}

\begin{proof}
The canonical functor from the right-hand side sends
$\widetilde x$ to the corresponding trivial torsor with
its standard identification over $k$.

Conversely, a deformation of $[x]$ consists of a torsor
$P\to\Spec A$, an equivariant map $P\to X$, and an
identification of its special fiber with the object
represented by $x$. This identification determines a
point $p_0\in P(k)$. Since $P\to\Spec A$ is smooth,
$p_0$ lifts to a point of $P(A)$, giving a trivialization
compatible with the prescribed special fiber.
The resulting map corresponds to an element of $X_x(A)$.
Arrows between two such objects are exactly the elements
of $G_1(A)$ carrying one lift to the other.
Thus, the canonical functor is fully faithful and
essentially surjective. Moreover, it is natural in $A$.
\end{proof}

Let $L$ be a finite-dimensional nilpotent dgLa.
Its \emph{Maurer--Cartan scheme} is defined by
$$\MC(L)(A)=\left\{\gamma\in L^1\otimes A:
 d\gamma+\frac12[\gamma,\gamma]=0\right\},$$
for every commutative $k$-algebra $A$.
The unipotent algebraic group $\exp(L^0)$ acts by gauge
transformations, defining the quotient stack
$$
 \mathfrak{MC}(L)=
 [\MC(L)/\exp(L^0)].
$$
For $\alpha\in\MC(L)(k)$, let
$$L^\alpha=(L,d^\alpha,[\cdot,\cdot]),\quad d^\alpha=d+[\alpha,-],
$$
be the corresponding twisted dgLa.

\begin{thm}\label{dglaconmc}
For every $\alpha\in\MC(L)(k)$, there is a natural
equivalence
$$
 \mathfrak{Def}_{\mathfrak{MC}(L),[\alpha]}
 \simeq
 \mathfrak{Def}_{L^\alpha}.
$$
Thus $L^\alpha$ controls the formal deformations of
$\mathfrak{MC}(L)$ at $[\alpha]$.
\end{thm}

\begin{proof}
Every lift of $\alpha$ over $A$ has the form
$\alpha+\eta$, with $\eta\in L^1\otimes\mathfrak m_A$.
Since $\alpha$ satisfies the Maurer--Cartan equation, we have
$$
 d(\alpha+\eta)+\frac12[\alpha+\eta,\alpha+\eta]
 =
 d^\alpha\eta+\frac12[\eta,\eta].
$$
Translation by $\alpha$ therefore gives a bijection
$$
 \MC(L^\alpha\otimes\mathfrak m_A)
 \xrightarrow\
 \MC(L)_\alpha(A).
$$
Moreover, we have
$$
 \ker\bigl(\exp(L^0)(A)\to\exp(L^0)(k)\bigr)
 =
 \exp(L^0\otimes\mathfrak m_A).
$$
Notice that for
$a\in L^0\otimes\mathfrak m_A$,
$$
 d^{\alpha}a=da-[a,\alpha],\quad e^a\star_{d}(\alpha+\eta)=\alpha+e^a\star_{d^\alpha}\eta,
$$
which means that translation intertwines the gauge actions.
\Cref{localquo} now yields
 \begin{align*}
 \mathfrak{Def}_{\mathfrak{MC}(L),[\alpha]}(A)
 &\simeq
 \MC(L)_\alpha(A)/\!/\exp(L^0\otimes\mathfrak m_A)\\
 &\cong
 \Del(L^\alpha\otimes\mathfrak m_A)
 =
 \mathfrak{Def}_{L^\alpha}(A).
 \end{align*}
Furthermore, the above proof is functorial in $A$.
\end{proof}

The following example indicates the role of the above machinery of dgLa and Maurer--Cartan formulations in the classical deformation problem of Lie algebra structures.

\begin{ex}
Let $V$ be an $n$-dimensional vector space with basis
$\{e_1,\ldots,e_n\}$, and let $\nu:\wedge^2 V\to V$ be a linear map. For $1\leq i< j\leq n$, write
$$
 \nu(e_i,e_j)=\sum_{a=1}^n c_{ij}^a e_a,
$$
where $\{c_{ij}^a\}$ denote the structure constants of $\nu$. The scheme of Lie brackets on $V$ is
$$
 \Lie(V)
 =
 \Spec\left(
 k[c_{ij}^a\mid 1\leq i<j\leq n]/(J_{ijk}^{\ell})
 \right),
$$
where, for $1\leq l\leq n$, $1\leq i< j<k\leq n$,
$$
 J_{ijk}^{\ell}
 =
 \sum_{a=1}^n
 \left(
 c_{ij}^a c_{ak}^{\ell}
 +c_{jk}^a c_{ai}^{\ell}
 +c_{ki}^a c_{aj}^{\ell}
 \right).
$$
For every commutative $k$-algebra $A$,
$$
 \Lie(V)(A)
 =
 \{\text{$A$-linear Lie brackets on }V\otimes A\}.
$$

The group $\GL(V)$ acts by transport of structure, and
$$
 \mathfrak{Lie}_n=[\Lie(V)/\GL(V)]
$$
is the moduli stack of Lie algebra structures on vector
bundles of rank $n$. For $A\in\Art$, every such vector
bundle is free, so $\mathfrak{Lie}_n(A)$ is equivalent
to the groupoid of Lie brackets on $V\otimes A$ and
$A$-linear Lie algebra isomorphisms.

Fix $\nu\in\Lie(V)(k)$. By Proposition~\ref{localquo},
$\mathfrak{Def}_{\mathfrak{Lie}_n,[\nu]}(A)$ is equivalent
to the groupoid whose objects are brackets $\nu_A$
satisfying
$$
 \nu_A\bmod\mathfrak m_A=\nu,
$$
and whose arrows $\nu_A\to\nu_A'$ are elements
$$
 g\in\GL(V\otimes A),
 \quad
 g\bmod\mathfrak m_A=\operatorname{Id}_V,
 \quad
 g\cdot\nu_A=\nu_A'.
$$

Let $L_{\mathrm{NR}}$ be the graded Lie algebra
$$
 L_{\mathrm{NR},V}^q
 =
 \Hom_k(\wedge^{q+1}V,V),
 \quad q\geq0,
$$
equipped with the Nijenhuis--Richardson bracket in
\Cref{NR} and zero differential.
The Jacobi identity is equivalent to
$[\nu,\nu]_{\mathrm{NR}}=0$, so $\nu$ is a
Maurer--Cartan element. Its twist is
$$
 L_{\mathrm{NR},V}^\nu
 =
 \bigl(L_{\mathrm{NR},V},[\nu,-]_{\mathrm{NR}},
 [\cdot,\cdot]_{\mathrm{NR}}\bigr).
$$
Writing $\nu_A=\nu+\eta$ identifies deformed brackets
with $\MC(L_{\mathrm{NR},V}^\nu\otimes\mathfrak m_A)$.
Furthermore, exponential and logarithm identify
$$
 \ker\bigl(\GL(V\otimes A)\to\GL(V)\bigr)
 =
 \exp\bigl(\End(V)\otimes\mathfrak m_A\bigr),
$$
and transport of structure agrees with the gauge action.
Consequently,
$$
 \mathfrak{Def}_{\mathfrak{Lie}_n,[\nu]}
 \simeq
 \mathfrak{Def}_{L_{\mathrm{NR},V}^\nu}.
$$
Thus the twisted Nijenhuis-Richardson dgLa controls
the formal deformations of the Lie algebra $(V,\nu)$.
\end{ex}

\section{Global extensions and non-abelian extensions of Lie algebras}\label{recall}
In this section, we recall the classification of global extensions of Lie algebras by global cohomology. Moreover, we will describe the classification result from the point of view of action groupoids. 

 Let $\g=(V,\mu)$ be a Lie algebra, let $V''$ be a vector space, and let $p:V''\to V$ be a surjective linear map. Denote by $V'=\ker p$.
 \begin{defi}{\rm (\cite[pp. 3--4]{AM15})}
\item A \emph{global extension structure} on $V''$ is a Lie bracket $\nu''$ for which
$$p:(V'',\nu'')\longrightarrow (V,\mu)$$
is a Lie algebra homomorphism. The set of global extension structures on $V''$ with respect to $p$ and $\g=(V,\mu)$ is denoted by $\mathrm{Gext}(\g, V')$. 

Two such structures $\nu''_1$ and $\nu''_2$ are \emph{equivalent (framed equivalent)}, if there is a Lie algebra isomorphism $\theta:(V'',\nu''_1)\to(V'',\nu''_2)$ such that
$$
\theta|_{V'}=\mathrm{id}_{V'},
\quad p\circ\theta=p.
$$
The set of equivalence classes is denoted by $\mathrm{GExt}(\g,V')$.
\end{defi}
However, in various interactions with deformation theory, the first condition is so strong, which will lead to additional rigidity results. To avoid this, we introduce a weaker notion of an equivalence relation, which we call the \emph{unframed equivalence}.

\begin{defi}
Two such structures $\nu''_1$ and $\nu''_2$ are \emph{unframed equivalent} if only the condition $p\circ\theta=p$ is imposed for the Lie algebra isomorphism $\theta$. The set of unframed equivalence classes is denoted by $\mathrm{GExt}^{\mathrm{un}}(\g,V')$.
\end{defi}
\begin{rmk}
A direct consequence of $p\circ\theta=p$ is that $\theta(V')=V'$ and $\theta|_{V'}$ is automatically a Lie algebra isomorphism from $(V',\nu''_{1}|_{V'})$ to the one induced by $(V',\nu''_{2}|_{V'})$. 
\end{rmk}

Throughout, a triple $(V'',p,\mathfrak g)$, where $\mathfrak g=(V,\mu)$ is a Lie algebra, $V''$ is a vector space, and $p\colon V''\rightarrow V$ is a surjective linear map, will be called a \emph{global extension datum}. We write $V'=\ker(p)$.

Let $m=\dim V''$ and $n=\dim\g$, and choose bases
$\{e_1,\ldots,e_m\}$ of $V''$ and $\{x_1,\ldots,x_n\}$ of $\g$.
Then $\mathrm{Gext}(\g,V')$ is an affine scheme with coordinate ring
$$
 \huaO_{\mathrm{Gext}(\g,V')}
 =k[t_1,\ldots,t_N]/\langle J,G\rangle,
 \quad N=\frac{m^2(m-1)}{2},
$$
where $J$ and $G$ encode the Jacobi identity and the condition that
$p$ is a Lie algebra homomorphism, respectively.

Consider the algebraic group
\begin{equation*}\label{eq:Gun-def}
 \mathrm{G}^{\mathrm{un}}_p=
 \left\{\theta\in\mathrm{GL}(V'')
 \;\middle|\;p\circ\theta=p
 \right\},
\end{equation*}
and its closed subgroup
\begin{equation*}\label{eq:G-def}
 \mathrm{G}_p=\left\{\theta\in\mathrm{GL}(V'')
 \;\middle|\;
 \theta|_{V'}=\operatorname{Id}_{V'},\quad p\circ\theta=p
 \right\},
\end{equation*}
which is the pointwise stabilizer of $V'$ in $G^\un_p$.

To describe these groups more explicitly, we fix a linear section $s:V\to V''$ of $p$.
\begin{pro}\label{prop:G-Hom}
The group $\mathrm{G}^{\mathrm{un}}_p$ is a closed algebraic subgroup of $\mathrm{GL}(V'')$. Moreover, the map
\begin{equation*}\label{eq:Phi-r}
 \Phi^{\un}: \Hom(V,V')\rtimes\GL(V')\longrightarrow\mathrm{G}^{\un}_p,\quad
 (r,B)\longmapsto\Phi^{\un}_{ (r,B)},
\end{equation*}
defined by
\begin{equation}\label{eq:Phi-r-formula}
 \Phi^{\un}_{(r, B)}(u+s(x))=B(u)+s(x)+r(x),\quad x\in V,\quad u\in V',
\end{equation}
is an isomorphism of algebraic groups. Here $\Hom(V,V')$ is
viewed as an additive vector group and the multiplication of $\Hom(V,V')\rtimes\GL(V')$ is given by 
$$(r,B)\cdot(r',B')=(r+Br', BB').$$
\end{pro}
\begin{proof}
Under the identification $V''\cong V'\oplus V$ induced by $s$,
we have the following realization 
$$
 \mathrm{G}^{\un}_p=\left\{\begin{pmatrix}
  B & r\\
  0 & I_V
 \end{pmatrix}\;\middle|\; r\in\Hom(V,V')\right\}.$$
Thus, $\mathrm{G}^{\un}_p$ is a closed subgroup of $\GL(V'')$. Moreover,
$\Phi^{\un}_{(r, B)}\Phi^{\un}_{(r', B')}=\Phi^{\un}_{(r+Br', BB')}$, and both $\Phi^{\un}$ and its inverse,
given by extracting the upper blocks, are regular.
Hence, $\Phi^{\un}$ is an isomorphism of algebraic groups.

Notice that $\mathrm{G}_p$ is the pointwise stabilizer of $V'$ in $\mathrm{G}_p^{\un}$, i.e., 
$$\mathrm{G}_p=\left\{\begin{pmatrix}
  I_{V’} & r\\
  0 & I_{V}
 \end{pmatrix}\;\middle|\; r\in\Hom(V,V')\right\},$$
 which implies that $\mathrm{G}_p$ is a unipotent group and $(\Phi^{\un})^{-1}(\mathrm{G}_p)=\Hom(V, V')\times\{I_{V'}\}$. Indeed, the map $\Phi^{\un}|_{\Hom(V, V')\times\{I_{V'}\}}$ does not depend on the choice of $s$.
\end{proof}

For $\theta\in\mathrm{G}^{\mathrm{un}}_p$ and
$\nu''\in\mathrm{Gext}(\g,V')$, define the transported bracket by
\begin{equation}\label{eq:G-action-X}
  (\theta\cdot\nu'')(a,b)
  = \theta\bigl(\nu''(\theta^{-1}a,\theta^{-1}b)\bigr),
  \quad \forall a, b\in V''.
\end{equation}
Since $p\circ\theta=p=p\circ\theta^{-1}$, the map $p$ remains a
homomorphism to $\g$ with respect to the transported bracket. Thus,
transport of structure defines an action of $\mathrm{G}^{\mathrm{un}}_p$ on $\mathrm{Gext}(\g,V')$ and, by restriction, an action of its subgroup $\mathrm{G}_p$. Denote the corresponding action groupoid by
$\mathrm{Gext}(\g,V')/\!/\mathrm{G}_p$ and $\mathrm{Gext}(\g,V')/\!/\mathrm{G}^{\un}_p$ respectively.
The sets $\mathrm{GExt}(\g,V')$ (resp. $\mathrm{GExt}^{\un}(\g,V')$) classify global extensions up to equivalence (resp. unframed equivalence). Hence, we obtain natural bijections
\begin{align*}
\mathrm{GExt}(\g,V')\cong
 \pi_0\big(\mathrm{Gext}(\g,V')/\!/\mathrm{G}_p\big),\quad
 \mathrm{GExt}^{\un}(\g,V')\cong
\pi_0\big(\mathrm{Gext}(\g,V')/\!/\mathrm{G}^{\un}_p\big).
\end{align*}

Next, we recall the concept of global cohomology of Lie algebras and use it to classify isomorphism classes of global extensions of Lie algebras. For more detail, we refer the reader to \cite{AM14-2, AM15}. 

\begin{defi}
Let $\g=(V, \mu)$ be a Lie algebra and $V'$ a vector space. Let $\pi_{V'}\in\Hom(\wedge^{2}V', V'),\ \varrho\in\Hom(V\otimes V', V'), $ and $\omega\in\Hom(\wedge^{2}V, V') $. The triple $(\pi_{V'},\varrho,\omega)$ is said to be a \emph{global 2-cocycle} of $\g$ with values in $V'$ if the following equations 
\begin{align}
\label{gc1} \sum_{\mathrm{cyc}}\pi_{V'}(\pi_{V'}(u,v),w)&=0,  \\
\label{gc4} \varrho(x)\pi_{V'}(u,v)&=\pi_{V'}(\varrho(x)u,v)+\pi_{V'}(u,\varrho(x)v),\\
\label{gc5} [\varrho(x),\varrho(y)]u-\varrho(\mu(x,y))u&=\pi_{V'}(\omega(x,y),u),\\
\label{gc6} \sum_{\mathrm{cyc}}\varrho(x)\omega(y,z)
&=\sum_{\mathrm{cyc}}\omega(\mu(x,y),z),
\end{align}
hold for all $x, y, z\in V$ and $u, v, w\in V'$. The set of these triples is denoted by $\mathrm{GZ}^2(\g,V')$.
\end{defi}
Notice that \eqref{gc1} shows that $\pi_{V'}$ defines a Lie algebra structure on $V'$, and \eqref{gc4} shows that $\varrho(x)$ is a derivation of the Lie algebra $(V', \pi_{V'})$ for all $x\in V$.

\begin{defi}
Two global 2-cocycles $(\pi_{V'},\varrho,\omega)$ and $(\pi_{V'}',
\varrho',\omega')$ are said to be \emph{ equivalent (framed equivalent)} if there exists a linear map $\tau:V \lon V'$ such that
\begin{align}
\label{eqg1}\pi_{V'}&=\pi_{V'}',\\
\label{eqg2}\varrho(x)u-\varrho'(x)u&=\pi_{V'}(\tau(x), u),\\
\label{eqg3}\omega(x, y)-\omega'(x, y)&=\varrho'(x)\tau(y)-\varrho'(y)\tau(x)+\pi_{V'}(\tau(x), \tau(y))-\tau(\mu(x, y)),
\end{align}
for all $x,y\in V$ and $u\in V'$. The set of equivalence classes is denoted by $\mathrm{GH}^{2}(\g, V')$, called the \emph{second global cohomology} of $\mathfrak{g}$ with values in $V'$.
\end{defi}

The set $\mathrm{GZ}^2(\g, V')$ of global $2$-cocycles is a closed
affine subscheme of
$$
 \Hom(\wedge^2V',V')
 \oplus\Hom(V\otimes V',V')
 \oplus\Hom(\wedge^2V,V').
$$
The additive group $\Hom(V,V')$ acts on $\mathrm{GZ}^2(\g,V')$ by
$$
 \tau\cdot(\pi_{V'},\varrho,\omega)=(\pi_{V'},\varrho^\tau,\omega^\tau),
$$
where
$$
 \begin{aligned}
 \varrho^\tau(x)u
 &=\varrho(x)u-\pi_{V'}(\tau(x),u),\\
 \omega^\tau(x,y)
 &=\omega(x,y)-\varrho(x)\tau(y)+\varrho(y)\tau(x)
   +\pi_{V'}(\tau(x),\tau(y))+\tau(\mu(x,y)).
 \end{aligned}
$$
Its orbits are precisely the equivalence classes of global
$2$-cocycles. Thus,
$$
 \mathrm{GH}^2(\g,V')\cong
 \pi_0\big(\mathrm{GZ}^2(\g,V')/\!/\Hom(V,V')\big),
$$
where $\mathrm{GZ}^2(\g,V')/\!/\Hom(V,V')$ denotes the action groupoid.

We now introduce a notion of unframed equivalence for global $2$-cocycles, analogous to that for global extensions.

\begin{defi}
Two global $2$-cocycles $(\pi_{V'},\varrho,\omega)$ and $(\pi_{V'}',
\varrho',\omega')$ are said to be \emph{unframed equivalent} if there exists a linear map $\tau:V \lon V'$ and $B\in \GL(V')$ such that
\begin{align*}
\pi_{V'}'(u, v)&=B\pi_{V'}(B^{-1}u, B^{-1}v),\\
\varrho'(x)u&=B\varrho(x)B^{-1}(u)-\pi'_{V'}(\tau(x), u),\\
\omega'(x, y)&=B\omega(x,y)-B\varrho(x)B^{-1}\tau(y)+B\varrho(y)B^{-1}\tau(x)+\pi'_{V'}(\tau(x),\tau(y))+\tau(\mu(x,y)),
\end{align*}
for all $x,y\in V$ and $u\in V'$. The set of unframed equivalence classes is denoted by $\mathrm{GH}^{\un,2}(\g, V')$, called the \emph{unframed second global cohomology} of $\mathfrak{g}$ with values in $V'$.
\end{defi}

Denote by $H_{V,V'}:=\Hom(V,V')\rtimes\GL(V')$ with the group law
\begin{equation*}
(\tau, B)(\eta, C)=(\tau+B\eta, BC),\quad \forall B,C\in \GL(V'), \tau,\eta\in \Hom(V,V').
\end{equation*}

There is a natural action of $H_{V,V'}$ on the set $\mathrm{GZ}^2(\g,V')$ by
\begin{equation*}
    (\tau, B)\cdot(\pi,\varrho,\omega)=(\pi',\varrho',\omega'), \quad \forall (\tau, B)\in H_{V,V'}
\end{equation*}
such that for $x,y\in V,\ u,v\in V'$,
\begin{align*}
    \pi'_{V'}(u,v)&=B\pi_{V'}(B^{-1}u,B^{-1}v),\\
    \varrho'(x)u&=B\varrho(x)B^{-1}(u)-\pi'_{V'}(\tau(x),u),\\
    \omega'(x,y)&=B\omega(x,y)-B\varrho(x)B^{-1}\tau(y)+B\varrho(y)B^{-1}\tau(x)+\pi'_{V'}(\tau(x),\tau(y))+\tau(\mu(x,y)).
\end{align*}
The orbits are precisely the unframed equivalence classes of global
$2$-cocycles, i.e.,
$$
 \mathrm{GH}^{\un,2}(\g,V')\cong
 \pi_0\big(\mathrm{GZ}^2(\g,V')/\!/H_{V,V'}\big),
$$
where $\mathrm{GZ}^2(\g,V')/\!/H_{V,V'}$ denotes the action groupoid.

Now fix the global extension datum $(V'',p,\g)$, write $V'=\ker p$,
and choose a linear section $s:V\to V''$ of $p$.
The section induces an isomorphism of affine schemes
\begin{equation}\label{defipsi}
 \Psi_s:\mathrm{Gext}(\g,V')
 \longrightarrow\mathrm{GZ}^2(\g,V'),
 \quad
 \nu''\longmapsto(\pi_{V'},\varrho,\omega),
\end{equation}
defined by
\begin{equation}\label{2co}
 \begin{aligned}
 \pi_{V'}(u,v)&=\nu''(u,v),\\
 \varrho(x)u&=\nu''(s(x),u),\\
 \omega(x,y)&=\nu''(s(x),s(y))-s(\mu(x,y)).
 \end{aligned}
\end{equation}
Its inverse assigns to $(\pi_{V'},\varrho,\omega)$ the bracket
on $V''\cong V'\oplus V$ given by
\begin{equation*}\label{eq:crossed-bracket}
 [(u, x),(v, y)]= \bigl(
 \pi_{V'}(u,v)+\varrho(x)v-\varrho(y)u+\omega(x,y),\ \mu(x,y)\bigr).
\end{equation*}

\begin{pro}\label{isogpd}
Let $s\colon V\to V''$ be a linear section of $p$, and let
$H_{V,V'}=\Hom(V,V')\rtimes\mathrm{GL}(V')$. The pair
$\bigl((\Phi^{\un})^{-1},\Psi_s\bigr)$ induces an isomorphism of
algebraic action groupoids
$$
  \mathrm{Gext}(\g,V')/\!/\mathrm{G}^{\un}_p
  \rightarrow
  \mathrm{GZ}^2(\g,V')/\!/H_{V,V'}.
$$
It sends an object $\nu''$ to $\Psi_s(\nu'')$ and an arrow
$(\theta,\nu'')$ to
$\bigl((\Phi^{\un})^{-1}(\theta),\Psi_s(\nu'')\bigr)$.
Restricting to the action of $\mathrm{G}_p$ gives an isomorphism
$$
  \mathrm{Gext}(\g,V')/\!/\mathrm{G}_p
  \rightarrow
  \mathrm{GZ}^2(\g,V')/\!/\Hom(V,V').
$$
\end{pro}

\begin{proof}
Since $\Phi^{\un}$ and $\Psi_s$ are algebraic isomorphisms, it suffices
to verify equivariance. Let $h=(r, B)\in H_{V,V'}$ and
$\nu''\in\mathrm{Gext}(\g,V')$. Write
$$
  \Psi_s(\nu'')=(\pi_{V'},\varrho,\omega),
  \quad
  \Psi_s\bigl(\Phi_h^{\un}\cdot\nu''\bigr)
  =(\pi_{V'}^h,\varrho^h,\omega^h).
$$
By \eqref{eq:Phi-r-formula}, the inverse of $\Phi_h^{\un}$ is given by
$$
  (\Phi_h^{\un})^{-1}\bigl(u+s(x)\bigr)
  =B^{-1}\bigl(u-r(x)\bigr)+s(x).
$$
Substituting this expression into the transport formula
\eqref{eq:G-action-X} gives, for $x,y\in V$ and $u,v\in V'$,
\begin{align*}
  \pi_{V'}^h(u,v)
    &= B\bigl(\pi_{V'}(B^{-1}u,B^{-1}v)\bigr), \\
  \varrho^h(x)u
    &= \bigl(B\varrho(x)B^{-1}\bigr)u
       -\pi_{V'}^h\bigl(r(x),u\bigr), \\
  \omega^h(x,y)
    &= B\omega(x,y)+r\bigl(\mu(x,y)\bigr)
       -\bigl(B\varrho(x)B^{-1}\bigr)r(y)
       +\bigl(B\varrho(y)B^{-1}\bigr)r(x)
       +\pi_{V'}^h\bigl(r(x),r(y)\bigr).
\end{align*}
These are precisely the components of $h\cdot(\pi_{V'},\varrho,\omega)$.
Hence
\begin{equation}\label{eq:equivariance}
  \Psi_s\bigl(\Phi_h^{\un}\cdot\nu''\bigr)
  =h\cdot\Psi_s(\nu'').
\end{equation}
This proves equivariance and yields the first groupoid isomorphism,
whose inverse is induced by $\bigl(\Phi^{\un},\Psi_s^{-1}\bigr)$.

Finally, $\Phi^{\un}$ identifies the subgroup
$\{(r, I_{V'})\mid r\in\Hom(V,V')\}$ with $\mathrm{G}_p$.
Restricting \eqref{eq:equivariance} to this subgroup and identifying
$(r, I_{V'})$ with $r$ gives the second isomorphism.
\end{proof}

Passing to connected components yields the following classifications.

\begin{thm}\label{iso-cohomology}
Let $(V'',p,\g)$ be a global extension datum with $V'=\ker p$.
\begin{enumerate}
\item 
There is a bijection
\begin{equation*}
    \mathrm{GExt}(\g,V')\cong\mathrm{GH}^2(\g,V').
\end{equation*}
\item Moreover, there is another bijection
\begin{equation*}
    \mathrm{GExt}^{\un}(\g,V')\cong\mathrm{GH}^{\un,2}(\g,V').
\end{equation*}
\end{enumerate}
\end{thm}

\begin{rmk}
The first part of the above theorem recovers Theorem 2.9 in \cite{AM15}. Furthermore, the above theorem holds regardless of the choice of the section $s:V\to V''$.
\end{rmk}
Indeed, any global extension structure $\nu''$ on $V''$ induces a non-abelian extension of $\g= (V,\mu)$ by $(\ker p,\nu''|_{\ker p})$, i.e.,
\[\begin{CD}
0@>>>(\ker p,\nu''|_{\ker p})@>i>>(V'',\nu'')@>p>>\g@>>>0.
\end{CD}\]

There is an obvious relation between global extensions and non-abelian extensions of Lie algebras. Indeed, such relation also has a cohomological shadow.

Let $\g=(V, \mu)$ and $\g'=(V', \mu_{V'})$ be Lie algebras. Let $0\lon\g'\overset{\iota}{\lon}\g''\overset{j}{\lon}\g\lon0$ be a non-abelian extension of $\g$ by $\g'$. 
\begin{defi}\label{defi111}{\rm (\cite[Definition 1.10]{Fre14})}
A \emph{non-abelian 2-cocycle} on $\g$ with values in $\g'$ is a pair $(\chi, \psi)$ of linear maps $\chi:\wedge^2V\lon V'$ and $\psi:V\lon\Der(V',\mu_{V'})$ such that for all $x,y\in V$,
\begin{align*}
 [\psi(x),\psi(y)]-\psi(\mu(x,y))
 &=\operatorname{ad}_{\chi(x,y)},\\
 \sum_{\mathrm{cyc}}\psi(x)\chi(y,z)
 &=\sum_{\mathrm{cyc}}\chi(\mu(x,y),z).
\end{align*}
Two cocycles $(\chi,\psi)$ and $(\chi',\psi')$ are
equivalent if there exists $\beta\in\Hom(V,V')$ such that
 \begin{align*}
 \psi'(x)&=\psi(x)+\operatorname{ad}_{\beta(x)},\\
 \chi'(x,y)
 &=\chi(x,y)+\psi(x)\beta(y)-\psi(y)\beta(x)
   -\beta(\mu(x,y))+\mu_{V'}(\beta(x),\beta(y)).
 \end{align*}
The set of cocycles is denoted by
$\mathrm{Z}^2_{\mathrm{nab}}(V, \mu; V', \mu_{V'})$. Its quotient by the above equivalence relation is the \emph{second non-abelian cohomology set} $\mathrm{H}^2_{\mathrm{nab}}(V, \mu; V', \mu_{V'})$.
\end{defi}

The map $(\chi,\psi)\mapsto(\mu_{V'},\psi,\chi)$ identifies
non-abelian $2$-cocycles with global $2$-cocycles having
Lie bracket $\mu_{V'}$ on the kernel. Moreover, this map is equivariant with respect to gauge actions on both sides. Since the global action preserves the Lie bracket on the kernel, we obtain the following decomposition.

\begin{thm}\rm{(\cite[Corollary 2.10]{AM15})}\label{geandnabe}
Let $\operatorname{Lie}(V')(k)$ be the set of Lie algebra structures
on $V'$. Then
$$
 \mathrm{GH}^2(\g,V')
 \cong
 \bigsqcup_{\mu_{V'}\in\operatorname{Lie}(V')(k)}
 \mathrm{H}^2_{\mathrm{nab}}\big(\g,(V',\mu_{V'})\big).
$$
\end{thm}

\begin{rmk}
Notice that in \cite{AM15}, the authors considered the global extension problem of non-commutative Poisson algebras. Here we adopt their results by taking the associative products of Poisson algebras and the corresponding global $2$-cocycle terms to be $0$.
\end{rmk}

\section{Global extensions of Lie algebras and Deligne groupoids}\label{4}
In this section, we construct a dgLa whose Maurer--Cartan scheme parametrizes
global extension structures and whose gauge action recovers
their equivalence relation. This identifies the corresponding
Deligne groupoid with the classification groupoid of global
extensions. Passing to quotient stacks, we show that the
formal deformations of a fixed extension are controlled by
the associated twisted dgLa, independently of the choice
of section.

We first recall a construction of gLas that will be used below. Let $U$ and $U'$ be vector spaces, and let
$(U,U')^{k,l}=\wedge^k U\otimes\wedge^l U'$.
For $W=U$ or $U'$, let $\iota_W:W\hookrightarrow U\oplus U'$
be the canonical inclusion. The \emph{lift} of a map
$f\in\Hom((U,U')^{k,l},W)$ is the element
$f^\dagger\in\Hom(\wedge^{k+l}(U\oplus U'), U\oplus U')$ defined by
$$
    f^\dagger(u_1+u'_1,\ldots,u_{k+l}+u'_{k+l})
    =
    \sum_{\sigma\in S(k,l)}(-1)^\sigma
    \iota_W\!\left(
        f(u_{\sigma(1)},\ldots,u_{\sigma(k)},
          u'_{\sigma(k+1)},\ldots,u'_{\sigma(k+l)})
    \right),
$$
where $S(k,l)$ denotes the set of $(k,l)$-shuffles and
$(-1)^\sigma$ is the signature of $\sigma$.

Together with the canonical decomposition
$
    \wedge^n(U\oplus U')
    \cong\oplus_{k+l=n}(U,U')^{k,l},
$
these lifts identify $\huaL_{NR,U\oplus U'}$ with
$$
    \oplus_{n\geq 1}\Big(\oplus_{k+l=n}
    \left(
        \Hom((U,U')^{k,l},U)
        \oplus
        \Hom((U,U')^{k,l},U')
    \right)\Big),
$$
where the summands with $k+l=n$ have degree $n-1$.
We equip this space with the graded Lie bracket transported
from the Nijenhuis--Richardson bracket under this identification.

\begin{lem}\label{lem1111}
With the above notations,
$$\bigg(\oplus^{+\infty}_{n=1}\Big((\oplus_{k+l=n}\Hom((U, U')^{k,l},U))\oplus(\oplus_{k+l=n}\Hom((U, U')^{k,l},U'))\Big), ~~[\cdot, \cdot]\bigg),$$ is a gLa, where $[\cdot,\cdot]$ is induced by the bracket $[\cdot,\cdot]_{\NR}$ on $\huaL_{NR,U\oplus U'}$ and satisfies
\begin{equation}\label{defibra}
[f, g]^{\dagger}=[f^{\dagger}, g^{\dagger}]_{\NR}.
\end{equation}
\end{lem}

\subsection{A dgLa for global extensions}
Let $(V'',p,\g)$ be a global extension datum and  $V'$ denote $\ker p$, where $\g=(V,\mu)$. Define two
graded vector spaces
$$
  \huaL_{V,V'}=\bigoplus_{n\in\mathbb Z}\huaL_{V,V'}^n,
  \quad
  \huaL_{V,V'}^{\un}
    =\bigoplus_{n\in\mathbb Z}\huaL_{V,V'}^{\un,n}.
$$
Both spaces vanish in negative degrees. In degree zero, set
$$
  \begin{aligned}
    \huaL_{V,V'}^0
      &= \Hom(V,V'), \\
    \huaL_{V,V'}^{\un,0}
      &=\Hom(V,V')\oplus\operatorname{End}(V').
  \end{aligned}
$$
For each $n\geq 1$, set
$$
  \huaL_{V,V'}^n
  = \huaL_{V,V'}^{\un,n}
  = \bigoplus_{\substack{i,j\geq 0\\ i+j=n+1}}
      \Hom\bigl(\wedge^i V\otimes\wedge^j V',\,V'\bigr).
$$
Define the Lie bracket and differential on $\huaL_{V,V'}$ and $\huaL_{V,V'}^{\un}$ as follows
\begin{align*}
    [\cdot, \cdot]_{\huaL_{V, V'}}: \huaL_{V, V'}\times\huaL_{V, V'}&\longrightarrow\huaL_{V, V'}, & d:\huaL_{V, V'}&\longrightarrow\huaL_{V, V'},\\
    (f,g)&\longmapsto[f, g]_{\huaL_{V, V'}}=[f, g],& f&\longmapsto d(f)=[\mu, f],
\end{align*}
for any $f\in\huaL_{V, V'}^k, g\in\huaL_{V, V'}^l$, where $[\cdot, \cdot]$ is given in Lemma \ref{lem1111}.
\begin{pro}\label{nondgla2}
With these operations, $\huaL_{V,V'}$ and $\huaL_{V,V'}^{\un}$ are dgLas.
\end{pro}

\begin{proof}
Since $\mu$ is a Maurer--Cartan element of $L_{\NR,V\oplus V'}$,
by \Cref{twistdgla}, there is a differential graded Lie algebra structure on $L_{\NR,V\oplus V'}$ with differential $d=[\mu,\cdot]$. By \eqref{defibra}, we conclude that $(L_{V, V'}, [\cdot, \cdot]_{L_{V, V'}}, d)$ and $(L_{V, V'}^{\un}, [\cdot, \cdot]_{L_{V, V'}}, d)$ are dgLas. 
\end{proof}

With these preparations in place, we can now characterize $\mathrm{GZ}^{2}(\g,V')$ in terms of Maurer--Cartan elements of $(\huaL_{V, V'}, [\cdot, \cdot]_{\huaL_{V, V'}}, d)$.

\begin{thm}\label{zmc1}
The set of Maurer--Cartan elements of $(\huaL_{V, V'}, [\cdot, \cdot]_{\huaL_{V, V'}}, d)$ is in bijection with $\mathrm{GZ}^{2}(\g,V')$, i.e.,
    \begin{equation*}
        \mathrm{GZ}^{2}(\g,V')\cong\mathrm{MC}(\huaL_{V, V'}),\quad (\pi_{V'},\varrho, \omega)\mapsto \omega+\varrho+\pi_{V'}, \quad \forall (\pi_{V'},\varrho, \omega)\in\mathrm{GZ}^{2}(\g,V').
    \end{equation*}
    \begin{proof}
       Suppose $\omega+\varrho+\pi_{V'}$ satisfies the Maurer--Cartan equation, for $\omega\in\Hom(\wedge^{2}V,V'),\ \varrho\in\Hom(V\otimes V',V')$ and  $\pi_{V'}\in\Hom(\wedge^{2}V',V')$. By \eqref{defibra}, we have
        $[\omega,\omega]=0, ~[\mu,\pi_{V'}]=0$. Then, by the expansion of the Maurer--Cartan equation for $\omega+\varrho+\pi_{V'}$ in $\huaL_{V, V'}$, we get   
        \begin{align}
            \label{gmc1}\Hom(\wedge^{3}V,V')&:[\mu,\omega]+[\omega,\varrho]=0,\\
            \label{gmc2}\Hom(\wedge^{2}V\otimes V',V')&:2( [\mu,\varrho]+[\omega,\pi_{V'}])+[\varrho,\varrho]=0,\\
            \label{gmc3}\Hom(V\otimes\wedge^{2}V',V')&:[\varrho,\pi_{V'}]=0,\\
            \label{gmc4}\Hom(\wedge^{3}V',V')&:[\pi_{V'},\pi_{V'}]=0.
        \end{align}
        Notice that \eqref{gmc1}--\eqref{gmc4} are equivalent to \eqref{gc6}--\eqref{gc1}. Therefore, the triple $(\pi_{V'}, \varrho,\omega)$ is a global 2-cocycle.

        Conversely, let $(\pi_{V'},\varrho,\omega)\in \mathrm{GZ}^{2}(\g,V')$ be a global 2-cocycle of $\mathfrak{g}$ with values in $V'$. Then define $\zeta=\omega+\varrho+\pi_{V'}\in \huaL_{V, V'}^{1}$. By direct calculation, one obtains that $\zeta$ is a Maurer--Cartan element of the dgLa $(\huaL_{V, V'}, [\cdot, \cdot]_{L_{V, V'}}, d)$.
    \end{proof}
\end{thm}

We can characterize the action groupoid $\mathrm{GZ}^{2}(\g,V')/\!/\Hom(V,V')$ introduced in \Cref{recall} in terms of the Deligne groupoid of the corresponding nilpotent dgLa, as follows.

\begin{thm}\label{imporantthm}
\begin{enumerate}
    \item The action groupoid of the global 2-cocycles is the same as the Deligne groupoid of the nilpotent dgLa $L_{V,V'}$, i.e.,
   \begin{equation*}
      \mathrm{GZ}^{2}(\g,V')/\!/\Hom(V,V')\cong \Del(L_{V,V'}).
   \end{equation*}
   \item Moreover, the set of connected components of $\mathrm{Del}(\huaL_{V, V'})$ is parameterized by the second global cohomology $\mathrm{GH}^2(\g, V')$, i.e.,
   \begin{equation*}
       \pi_0(\mathrm{Del}(\huaL_{V, V'}))\cong \mathrm{GH}^2(\g, V').
   \end{equation*}
\end{enumerate}
\end{thm}
    \begin{proof}
 By Theorem \ref{zmc1}, the set of objects of  $\mathrm{Del}(\huaL_{V, V'})$ is $\mathrm{GZ}^{2}(\g, V')$. Then, we only need to prove that
 $\Hom_{\mathrm{Del}(\huaL_{V, V'})}(\zeta, \zeta')$ is in one-to-one correspondence with the Hom set between $\zeta=\omega+\varrho+\pi_{V'}$ and $ \zeta'=\omega'+\varrho'+\pi_{V'}'\in \mathrm{GZ}^{2}(\g, V')$. 
    
Let $\tau\in \Hom_{\mathrm{Del}(\huaL_{V, V'})}(\zeta, \zeta')$, such that
        \begin{equation}\label{gmceq}
            \zeta'=e^{\ad\tau}\zeta+\frac{\Id-e^{\ad\tau}}{\ad\tau}d(\tau).
        \end{equation}
        By \eqref{defibra}, we obtain $[\tau,\omega]=0$ and $[\tau,[\tau,\varrho]]=0$, which implies that
        \begin{equation*}
    e^{\ad\tau}\zeta=\omega+\varrho+\pi_{V'}+[\tau,\varrho]+[\tau, \pi_{V'}]+\frac{1}{2}[\tau,[\tau,\pi_{V'}]].
        \end{equation*}
        By \eqref{defibra}, we also obtain $[\tau,[\mu,\tau]]=0$, 
        which implies that
        \begin{eqnarray*}
            \frac{\Id-e^{\ad\tau}}{\ad\tau}d(\tau)=-\sum_{n\geq0}\frac{1}{(n+1)!}\ad^{n}\tau(d(\tau))=-\Big(d(\tau)+\frac{1}{2}[\tau,d(\tau)]+\cdots\Big)=-[\mu,\tau],
        \end{eqnarray*}
        where the summation is finite because $L_{V,V'}$ is nilpotent.
        We can rewrite equation (\ref{gmceq}) and rearrange the result according to the corresponding components. Then we obtain the following equations
        \begin{align*}
            \Hom(\wedge^{2}V,V')&:~~ \omega'=\omega+[\tau,\varrho]-[\mu,\tau]+\frac{1}{2}[\tau,[\tau,\pi_{V'}]],\\
            \Hom(V\otimes V',V')&:~~ \varrho'=\varrho+[\tau, \pi_{V'}],\\
            \Hom(\wedge^{2}V',V')&:~~\pi_{V'}'=\pi_{V'}.
        \end{align*}
        The above three equations give us the equalities (\ref{eqg3}), (\ref{eqg2}), and (\ref{eqg1}) respectively. Hence, we get 
        \begin{equation*}
            (\pi_{V'}',\varrho',\omega')=\tau\cdot(\pi_{V'},\varrho, \omega)
        \end{equation*} 

        Conversely, let $(\pi_{V'}',\varrho',\omega')=\tau\cdot(\pi_{V'},\varrho, \omega)$ in $\mathrm{GZ}^{2}(\g,V')$, i.e., the linear map $\tau:V\lon V'$ satisfies the equations (\ref{eqg1}), (\ref{eqg2}), and (\ref{eqg3}). Then by direct calculation, one can see that $\zeta\coloneqq \omega+\varrho+\pi_{V'}$ and $\zeta'\coloneqq \omega'+\varrho'+\pi_{V'}'$ are gauge equivalent through $e^{\ad\tau}$, which means that $\tau\in \Hom_{\mathrm{Del}(\huaL_{V, V'})}(\zeta, \zeta')$.

        By the fact that $L_{V,V'}^{0}$ is abelian, the Baker-Campbell-Hausdorff multiplication reduces to the addition of linear maps. We get 
        $$\exp(L_{V,V'}^{0})\cong \big(\Hom(V,V'),+\big)$$
        as algebraic groups.

        The equivariance follows from summarizing the above computations into the following equality:
        \begin{equation*}
            e^{\tau}\star(\omega+\varrho+\pi_{V'})=\tau\cdot(\pi_{V'},\varrho,\omega).
        \end{equation*}
        Therefore, we conclude that $\mathrm{GZ}^{2}(\g,V')/\!/\Hom(V,V')\cong \Del(L_{V,V'})$ as algebraic groupoids.
        
        The second statement follows from taking $\pi_0$ of the above groupoids.  
\end{proof}

\subsection{Formal deformations and independence of the section}
Introduce the quotient stacks
\begin{equation*}
 \mathcal Z_{\g,V'}:=[\mathrm{GZ}^2(\g,V')/\Hom(V, V')],
 \end{equation*}
 \begin{equation*}
 \mathcal E_p:=[\mathrm{Gext}(\g,V')/\mathrm{G}_p],
\end{equation*}
\begin{equation*}
\mathcal E_p^{\un}:=[\mathrm{Gext}(\g,V')/\mathrm{G}^{\un}_p].
\end{equation*}
Write
$$
 \mathfrak{MC}(L_{V, V'}):=[\MC(L_{V,V'})/\exp(L_{V, V'}^0)].
$$
For a stack $\mathfrak X$ and an object $x\in\mathfrak X(k)$,
we use the formal deformation functor in \Cref{fordef}
$$
 \mathfrak{Def}_{\mathfrak X,x}(A)=\mathfrak X(A)\times^h_{\mathfrak X(k)}\{x\},\quad A\in\Art.
$$
\begin{thm}\label{dgla2co}
There is an equivalence of stacks
$$
 \mathcal Z_{\g,V'}\simeq\mathfrak{MC}(L_{V, V'}).
$$
For a global cocycle $z=(\pi,\varrho,\omega)$, let
$\zeta=\pi+\varrho+\omega$. Then
$$
 \mathfrak{Def}_{\mathcal Z_{\g,V'},[z]} \simeq
 \mathfrak{Def}_{L_{V, V'}^\zeta},
$$
naturally as groupoid-valued functors on $\Art$. Therefore, the formal deformation problem of global 2-cocycles is controlled by the dgLa $L^{\zeta}_{V,V'}$.
\end{thm}
\begin{proof}
By \Cref{imporantthm}, there exists an equivalence of
quotient stacks, which carries $[z]$ to $[\zeta]$. Hence,
\Cref{dglaconmc} implies
$$
 \mathfrak{Def}_{\mathcal Z_{\g,V'},[z]}
 \simeq
 \mathfrak{Def}_{\mathfrak{MC}(L_{V, V'}),[\zeta]}
 \simeq
 \mathfrak{Def}_{L_{V, V'}^\zeta}.
$$
\end{proof}

\begin{cor}\label{defglobal2co}
A linear section $s:V\to V''$ induces an equivalence
$$
 \mathcal E_p\simeq\mathcal Z_{\g,V'}.
$$
For $\mu''\in\mathrm{Gext}(\g,V'')(k)$ and
$z_s=\Psi_s(\mu'')$, where $\Psi_s$ is given by \eqref{defipsi}, it induces a natural equivalence
$$
 \mathfrak{Def}_{\mathcal E_p,[\mu'']}
 \simeq
 \mathfrak{Def}_{\mathcal Z_{\g,V'},[z_s]}.
$$
\end{cor}
\begin{proof}
By \Cref{isogpd}, there is an isomorphism of
algebraic action groupoids carrying $\mu''$ to $z_s$.
Passing to quotient stacks and taking homotopy fibers
at the corresponding objects, we prove the result.
\end{proof}

Note that we need to choose a section $s$ in \Cref{defglobal2co}. However, the following result shows that it is independent of the choice when we study $\mathfrak{Def}_{\mathcal E_p, [\mu^{''}]}$ via dgLa.
 
 For a global extension structure $\mu''\in{\rm Gext}(\g, V'')$, take a section $s: V\lon V''$ of $p$. Denote by
 \begin{equation*}\label{twdgla}
     \huaL_{V,V'}^{s}\coloneqq \huaL_{V,V'}^{(\pi_{V'}^{s},~\varrho^{s},~\omega^{s})}
     \quad
     \huaL_{V,V'}^{\un,s}\coloneqq (\huaL_{V,V'}^{\un})^{(\pi_{V'}^{s},~\varrho^{s},~\omega^{s})}
 \end{equation*}
 the twisted dgLas of $\huaL_{V,V'}$ and $\huaL_{V,V'}^{\un}$ by the corresponding Maurer--Cartan element $\omega^{s}+\varrho^{s}+\pi_{V'}^{s}$ of the global 2-cocycle $(\pi_{V'}^{s},\varrho^{s},\omega^{s})$ respectively, where $(\pi_{V'}^{s},\varrho^{s},\omega^{s})$ is given by \eqref{2co}. Let $s':V\to V''$ be another linear section of $p$. Define $\tau:V\to V'$ by 
\begin{equation}\label{eqdef:Ttau}
   \tau(x)=s(x)-s'(x), \quad \forall x\in V.
\end{equation}
 
\begin{pro}\label{prop:section-independence}
The dgLas $\huaL_{V,V'}^{s}$ and $\huaL_{V,V'}^{s'}$ are isomorphic.
   \end{pro}
    \begin{proof}
       Since $s, s'$ are sections of $p$, it is shown in \cite[Lemma 2.6]{AM15} that the triples $(\pi_{V'}^{s},\varrho^{s},\omega^{s})$ and $(\pi_{V'}^{s'},\varrho^{s'},\omega^{s'})$ are equivalent global $2$-cocycles. As in \Cref{imporantthm}, $\omega^{s}+\varrho^{s}+\pi_{V'}^{s}$ and $\omega^{s'}+\varrho^{s'}+\pi_{V'}^{s'}$ are gauge equivalent Maurer--Cartan elements of $\huaL_{V,V'}$ via $e^{\tau}$. By \cite[Corollary 5.10]{KS}, the two dgLas are isomorphic.
    \end{proof}

\begin{thm}\label{dgladefglobal}
For every global extension structure $\mu''$ and every
linear section $s$ of $p$, there is a natural equivalence
$$
 \mathfrak{Def}_{\mathcal E_p,[\mu'']}
 \simeq
 \mathfrak{Def}_{L_{V, V'}^s}.
$$
Explicitly, for $A\in\Art$,
$$
 \mathfrak{Def}_{\mathcal E_p,[\mu'']}(A)
 \simeq
 \Del(L_{V, V'}^s\otimes\mathfrak m_A).
$$
The controlling dgLa is independent of $s$ up to the
isomorphisms in Proposition~\ref{prop:section-independence}.
\end{thm}
\begin{proof}
By \Cref{dgla2co} and \Cref{defglobal2co}, there are natural equivalences
 \begin{align*}
 \mathfrak{Def}_{\mathcal E_p,[\mu'']}
 \simeq
 \mathfrak{Def}_{\mathcal Z_{\g,V'},[z_s]}
 \simeq
 \mathfrak{Def}_{\mathfrak{MC}(L_{V, V'}),[\zeta_s]}
 \simeq
 \mathfrak{Def}_{L_{V, V'}^s}.
 \end{align*}
\end{proof}
Similarly, for the unframed equivalence, we also have the following characterization.

\begin{thm}\label{thm:unframed-stack}
The stack of global extension structures with unframed isomorphisms has the two isomorphic presentations
\begin{equation*}
\mathcal E_p^{\un}=[\mathrm{Gext}(\g,V')/\mathrm{G}^{\un}_p]\cong[\MC(L^{\un}_{V,V'})/H_{V,V'}].
\end{equation*}
Every $\mathrm{G}^{\un}_p$-torsor over $\operatorname{Spec}\K$ is trivial. Consequently, after choosing the section of $p$, there are
isomorphisms of groupoids
\begin{equation*}
    \mathcal E_p^{\un}(\K) \simeq\mathrm{Gext}(\g,V')(\K)/\!/ \mathrm{G}^{\un}_p(\K) \cong \MC(L^{\un}_{V,V'})(k)/\!/ H_{V,V'}(\K),
\end{equation*}
and the set of connected components is $\mathrm{GExt}^{\mathrm{un}}(\g,V')$.
The natural map sending a global extension to the Lie algebra structure on the kernel is equivariant and induces a morphism
\begin{equation}\label{eq:unframed-kernel-map}
\widetilde q:\mathcal E_p^{\un}
\longrightarrow[\Lie(V')/\GL(V')].
\end{equation}
\end{thm}

\begin{proof}
The intrinsic action of $\mathrm{G}^{\un}_p$ becomes the action of $H_{V,V'}$ under the chosen section, which proves the equivalence of quotient stacks. Consider the split exact sequence
\begin{equation*}
    1\longrightarrow \Hom(V,V')\longrightarrow H_{V,V'}
\longrightarrow\GL(V')\longrightarrow1.
\end{equation*}
By the associated long exact sequence of pointed sets, we have 
\begin{equation*}
    1\longrightarrow H^1(k,\Hom(V,V'))\longrightarrow H^1(k, H_{V,V'})
\longrightarrow H^1(k,\GL(V'))=1.
\end{equation*}
which implies that $H^1(k,\mathrm{G}^{\un}_p)$ is trivial, since $\Hom(V,V')$ is unipotent and has trivial $H^1$. Hence, every $\mathrm{G}^{\un}_p$-torsor over $\operatorname{Spec} k$ is trivial, which further implies that 
$$\mathcal E_p^{\un}(\K) \simeq\mathrm{Gext}(\g,V')(\K)/\!/ \mathrm{G}^{\un}_p(\K).$$

For $(\tau, B)\in H_{V,V'}$ and $\zeta\in \MC(L^{\un}_{V,V'})$, restricting the transported bracket to $V'$, we have
\begin{equation*}
    q((\tau,B)\cdot \zeta)=B\cdot q(\zeta),\quad
(B\cdot\pi)(u,v)=B\pi(B^{-1}u,B^{-1}v).
\end{equation*}
Thus, the map $q:\MC(L^{\un}_{V,V'})\to\Lie(V')$ is equivariant for the homomorphism $H_{V,V'}\to\GL(V')$, and it descends to \eqref{eq:unframed-kernel-map}.
\end{proof}

\begin{thm}\label{unframed-formal-moduli} $L_{V, V'}^{\un}$ controls the formal deformation of global extension structures under the unframed equivalence, i.e.,
    \begin{equation*}
        \mathfrak{Def}_{\mathcal E^{\un}_p,[\mu'']} \simeq
        \mathfrak{Def}_{L_{V, V'}^{\un,s}}.
    \end{equation*}
\end{thm}
\begin{proof}
According to \Cref{thm:unframed-stack} and \Cref{localquo}, we obtain 
\begin{equation*}
    \mathfrak{Def}_{\mathcal E^{\un}_p,[\mu'']}(A)\simeq \MC(L^{\un}_{V,V'})_{\mu''}(A)/\!/H_{1}(A),
\end{equation*}
where $H_{1}(A)=\ker\bigl(H_{V,V'}(A)\longrightarrow H_{V,V'}(\K)\bigr)$. Since $\MC(L^{\un}_{V,V'})=\MC(L_{V,V'})$, by \Cref{dgladefglobal}, we have
\begin{equation*}
    \MC(L^{\un}_{V,V'})_{\mu''}(A)=\MC(L_{V,V'})_{\mu''}(A)=\MC(L_{V,V'}^{s}\otimes\frkm_{A})=\MC(L_{V,V'}^{\un,s}\otimes\frkm_{A}).
\end{equation*} 
Since $\m_A$ is nilpotent and $\K$ has characteristic zero, the finite exponential and logarithm series identify $H_{1}(A)$ with $\exp(L_{V,V'}^{\un,0}\otimes\m_A)$. Thus, we obtain that  
\begin{equation*}\label{eq:unframed-formal-moduli}
\mathfrak{Def}_{\mathcal E^{\un}_p,[\mu'']}(A) \simeq
\Del\bigl(L_{V,V'}^{\un,s}\otimes\m_A\bigr),
\end{equation*}
which is natural in $A$.
\end{proof}

\subsection{Examples}
At the end of this section, we classify global extensions up to framed equivalence for two simple examples based on \Cref{zmc1} and \Cref{imporantthm}, which shows the effectiveness of our machinery.
\begin{ex}\label{ex11}
Let $\g=(V, \mu)$ be a finite-dimensional abelian Lie algebra over $k$. Denote by $k\oplus V\overset{p}{\lon}V$ the projection, i.e., 
$$p(r,x)=x, \quad \forall x\in V,r\in k.$$
In the following, we describe the set $\mathrm {GExt}(\mathfrak{g},k)$ of global extensions on $k\oplus V$ with respect to $p$ and $V$. For $\omega+\varrho+\pi_{k}\in\huaL^{1}_{V,k}$, where
$$\pi_k\in\Hom(\wedge^{2}k,k)=\{0\},\quad\varrho\in\Hom(V\otimes k, k)\cong V^*, \quad\omega\in\Hom(\wedge^{2}V, k)\cong\wedge^{2}V^*, $$
we have that $\omega+\varrho+\pi_k\in\mathrm{MC}(\huaL_{V, k})$ if and only if 
    $$
    \varrho(x)\omega(y, z)+\varrho(y)\omega(z, x)+\varrho(z)\omega(x, y)=0, \quad \forall x,y,z\in V.
    $$
That is, $\varrho\wedge\omega=0$.
\begin{itemize}
\item If $\varrho\neq 0$, there exists $\varsigma\in V^*$ such that $\omega=\varrho\wedge\varsigma$. \item If $\varrho=0$, every $\omega\in\Hom(\wedge^{2}V, k)$ satisfies $\varrho\wedge\omega=0$.  
\end{itemize}
Therefore,

\begin{align*}
    \mathrm{GZ}^{2}(V, k)\cong\mathrm{MC}(\huaL_{V, k})&\cong\Big(\bigsqcup_{\varrho\in\{ V^{*}\setminus0\}}\{\varrho\}\times (V^{*}/k\varrho)\Big)\bigsqcup \Big(\{0\}\times \wedge^{2}V^{*}\Big)\\
    &\cong\Big(\bigsqcup_{\varrho\in\{ V^{*}\setminus0\}}\{\varrho\}\times k^{\dim V-1}\Big)\bigsqcup \Big(\{0\}\times k^{\frac{(\dim V-1)\dim V}{2}}\Big).
\end{align*} 
    Moreover, by the definition of gauge equivalence of {\rm Maurer--Cartan} elements, $\omega+\varrho$ and $\omega'+\varrho'$ are equivalent if and only if 
    $$
    \varrho=\varrho'=0, \quad \omega=\omega',
    $$
or there exists a linear map $\tau\in V^*$ such that
    $$
     \varrho=\varrho'\neq0, \quad \omega-\omega'=\varrho'\wedge\tau.
    $$
     Thus, by \Cref{iso-cohomology}(1) and \Cref{imporantthm}, there are bijections
    $$
    \mathrm {GExt}(\mathfrak{g},k)\cong\pi_0({\rm Del}(\huaL_{V, k}))\cong\mathrm{GH}^{2}(\mathfrak{g}, k)\cong(V^*\setminus0)\bigsqcup k^{\frac{(\dim V-1)\dim V}{2}}.
    $$
    Moreover, the stabilizer of $(\varrho,\omega)$ is $k\varrho$ if $\varrho\neq0$, and is $V^{*}$ if $\varrho=0$. 
\end{ex}

\begin{ex}\label{exlied}
Let $\g=(V, \mu)$ be a one-dimensional Lie algebra over $k$, and let $V'$ be a finite-dimensional vector space over $k$. Denote by $V'\oplus V\overset{p}{\lon}V$ the projection, i.e. 
$$
p(u,x)=x, \quad \forall x\in V,u\in V'.
$$
For the triple $\omega+\varrho+\pi_{V'}\in\huaL^{1}_{V, V'}$, where
    \begin{equation*}
       \pi_{V'}\in\Hom(\wedge^{2}V', V'),\quad\varrho\in\Hom(V\otimes V', V')\cong\Hom(V',V'), \quad\omega\in\Hom(\wedge^{2}V, V')=\{0\},
    \end{equation*}
   we have that $\varrho+\pi_{V'}\in{\rm MC}(\huaL_{V, V'})$ if and only if 
   \begin{align*}
       0&=\pi_{V'}(\pi_{V'}(u, v), w)+\pi_{V'}(\pi_{V'}(v, w), u)+\pi_{V'}(\pi_{V'}(w, u), v)\\
       \varrho(\pi_{V'}(u, v))&=\pi_{V'}(\varrho(u), v)+\pi_{V'}(u, \varrho(v)),
   \end{align*}
which means that $\pi_{V'}$ is a Lie bracket on $V'$ and $\varrho$ is a derivation of $\pi_{V'}$. Thus, 
$$
\mathrm{GZ}^{2}(\g, V')\cong{\rm MC}(\huaL_{V, V'})\cong\bigsqcup_{\pi_{V'}\in{\rm MC}(\huaL_{V'})}\Der(V',\pi_{V'}),
$$
where $\Der(V', \pi_{V'})$ is the set of derivations of the Lie algebra $(V', \pi_{V'})$.

    If $\varrho+\pi_{V'}$ and $\varrho'+\pi_{V'}'$ are equivalent in ${\rm MC}(\huaL_{V, V'})$, then there exists $\tau\in \Hom(V, V')\cong V'$ such that 
    \begin{equation*}
        \pi_{V'}=\pi_{V'}', \quad  \varrho=\varrho'+\ad (\tau).
    \end{equation*}
    Hence, by \Cref{iso-cohomology}(1), there is
    $$
    \mathrm {GExt}(\mathfrak{g}, V')\cong\mathrm{GH}^{2}(\g, V')=\bigsqcup_{\pi_{V'}\in{\rm MC}(\huaL_{V'})}\Out(V',\pi_{V'}),$$ where $\Out(V',\pi_{V'})$ is the vector space of outer derivations of the Lie algebra $(V',\pi_{V'})$. 
\end{ex}

\section{Non-abelian extensions of dgLas and fibrations of groupoids}\label{sec:non}
In this section, we investigate the relation between the Deligne groupoid corresponding to the global cohomology of Lie algebras and the one corresponding to the non-abelian cohomology of Lie algebras.

We first recall the relation between the Deligne groupoid of a specific dgLa and the second non-abelian cohomology of Lie algebras.

\begin{pro}\label{nondgla1}{\rm (\cite[Lemma 2.1, Corollary 2.1, Proposition 2.2]{Fre14})}
\begin{enumerate}
\item
Let $\g=(V, \mu)$ and $\g'=(V', \mu')$ be Lie algebras. Define a graded vector space $\widetilde{\huaL}_{V, V'}=\oplus_{n\in\mathbb{Z}}\widetilde\huaL_{V, V'}^n$ by 
\begin{equation*}
\widetilde\huaL_{V, V'}^n=
\begin{cases}
0, & n\leq-1;\\
\Hom(V, V'), & n=0;\\
\oplus_{\substack{i+j=n+1 \\ i\geq 1}}\Hom(\wedge^{i}V\otimes\wedge^{j}V', V'), & n\geq 1.
\end{cases}
\end{equation*}
Give the bracket and differential on $\widetilde\huaL_{V,V'}$ as follows
\begin{align*}
    [\cdot, \cdot]_{\widetilde\huaL_{V, V'}}: \widetilde\huaL_{V, V'}\times\widetilde\huaL_{V, V'}&\longrightarrow\widetilde\huaL_{V, V'}, &\dM:\widetilde\huaL_{V, V'}&\longrightarrow\widetilde\huaL_{V, V'},\\
    (f,g)&\longmapsto[f, g]_{\widetilde\huaL_{V, V'}}=[f, g],& f&\longmapsto \dM(f)=[\mu+\mu', f],
\end{align*}
for any $f\in\widetilde\huaL_{V, V'}^k, g\in\widetilde\huaL_{V, V'}^l$, where $[\cdot, \cdot]$ is given in Lemma \ref{lem1111}.
Then, the dgLa $(\widetilde\huaL_{V, V'}, [\cdot, \cdot]_{\widetilde\huaL_{V, V'}}, \dM)$ is nilpotent.
\item
The set of connected components of Deligne groupoid of $(\widetilde\huaL_{V, V'}, [\cdot, \cdot]_{\widetilde\huaL_{V, V'}}, \dM)$ is isomorphic to the second non-abelian cohomology of $\g=(V, \mu)$ with values in $\g'=(V', \mu')$, i.e., 
\begin{equation*}
    \pi_0\big(\mathrm{Del}(\widetilde\huaL_{V, V'})\big)\cong \mathrm{H}^2_{\mathrm{nab}}(V, \mu; V', \mu').
\end{equation*}
\end{enumerate}
\end{pro}

Indeed, the three dgLas $(\huaL_{V'}, [\cdot, \cdot]_{\NR}),\  (\widetilde{\huaL}_{V, V'}, [\cdot, \cdot]_{\widetilde{\huaL}_{V, V'}}, \dM)$ and $(\huaL_{V, V'}, [\cdot, \cdot]_{\huaL_{V, V'}}, d)$, introduced respectively in \Cref{dglamclie}, \Cref{nondgla1} and \Cref{nondgla2}, fit together into a short exact sequence of dgLas, as shown below. Here, in the construction of $(\widetilde{\huaL}_{V, V'}, [\cdot, \cdot]_{\widetilde{\huaL}_{V, V'}}, \dM)$, $V'$ is regarded as a trivial Lie algebra.

\begin{lem}\label{nonabdgla}
With the above notations, $(\huaL_{V, V'}, [\cdot, \cdot]_{\huaL_{V, V'}}, d)$ is a \emph{non-abelian extension} of $(\huaL_{V'}, [\cdot, \cdot]_{\NR})$ by $(\widetilde{\huaL}_{V, V'}, [\cdot, \cdot]_{\widetilde{\huaL}_{V, V'}}, \dM)$, i.e., the following sequence is exact
    \[\begin{CD}
    0@>>>\widetilde{\huaL}_{V,V'}@>\iota>>\huaL_{V,V'}@>\mathrm{pr}>>\huaL_{V'}@>>>0.
    \end{CD}\]
    where 
    \begin{equation*}
        \iota_{n}:\widetilde{\huaL}_{V,V'}^{n}\lon\huaL_{V,V'}^{n},\quad f\mapsto (f,0),\quad \forall f\in \oplus_{\substack{i+j=n+1 \\ i\geq 1}}\Hom(\wedge^{i}V\otimes\wedge^{j}V',V'),
    \end{equation*}
    and
    \begin{equation*}
        \mathrm{pr}_{n}:\huaL_{V,V'}^{n}\lon\huaL_{V'}^{n},\quad (g_{0}, \cdots g_{n+1})\mapsto g_{n+1},\quad \forall g_{i}\in \Hom(\wedge^{n+1-i}V\otimes\wedge^{i}V',V'),
    \end{equation*}
    with $n\geq1$. In the degree $0$ case, the inclusion is identity and the projection is $0$.
\end{lem}
By the fact that the three dgLas are nilpotent and Proposition \ref{indufunc}, the maps $$\iota:(\widetilde{\huaL}_{V, V'}, [\cdot, \cdot]_{\widetilde{\huaL}_{V, V'}}, \dM)\lon(\huaL_{V, V'}, [\cdot, \cdot]_{\huaL_{V, V'}}, d), \quad\   
\mathrm{pr}: (\huaL_{V, V'}, [\cdot, \cdot]_{\huaL_{V, V'}}, d)\lon(\huaL_{V'}, [\cdot, \cdot]_{\NR}),$$ 
induce morphisms of Deligne groupoids $$\iota_*:\Del(\widetilde{\huaL}_{V,V'})\lon\Del(\huaL_{V,V'}),\quad \mathrm{pr}_{*}:\Del(\huaL_{V,V'})\lon\Del(\huaL_{V'}).$$ At the level of groupoids, we have the following results.
\begin{pro}\label{pro:fiberation-deligne-gpd}
    (1) $\iota_{*}:\Del(\widetilde{\huaL}_{V,V'})\lon\Del(\huaL_{V,V'})$ is an inclusion of groupoids,

    (2) $\mathrm{pr}_{*}:\Del(\huaL_{V,V'})\lon\Del(\huaL_{V'})$ is a fibration of groupoids. 
    \begin{proof}
        By the fact that $\iota$ is injective on $\widetilde{\huaL}_{V,V'}^{0}$ and $\widetilde{\huaL}_{V,V'}^{1}$, the first statement is obvious.

       Notice that $\mathrm{Im} \mathrm{pr}_{*}\big(\mathrm{Obj}\big(\Del(\huaL_{V,V'})\big)\big)=\mathrm{Obj}(\Del(\huaL_{V'}))=\mathrm{MC}(\huaL_{V'})$. Suppose $\pi\in \mathrm{MC}(\huaL_{V'})$, then there exists an object $\kappa\in\mathrm{Obj}\big(\Del(\huaL_{V, V'})\big)$ such that $$\mathrm{pr}_{*}(\kappa)=\mathrm{pr}(\kappa)=\pi.$$ Let $\psi\in \St_{\pi}\big(\Del(\huaL_{V'})\big)$ be a morphism with source $\pi$ and target $\pi'$. Then there exists $a\in\huaL_{V'}^0$ such that $\psi=e^{a}$. By the surjectivity of $\mathrm{pr}$ on $\huaL_{V,V'}^{0}$, there exists an element $\tilde{a}\in \huaL_{V,V'}^{0}$ such that $\mathrm{pr}(\tilde{a})=a$. Define $\tilde{\psi}=e^{\tilde{a}}$, we have
        \begin{equation*}
            \mathrm{pr}_*(e^{\tilde{a}})=e^{\mathrm{pr}(\tilde{a})}=e^{a}.
        \end{equation*}
        Define $\kappa'=e^{\tilde{a}}\star \kappa$, we obtain that
        \begin{equation*}
            \mathrm{pr}_{*}(\kappa')=\mathrm{pr}(e^{\tilde{a}}\star \kappa)=e^{\mathrm{pr}(\tilde{a})}\star \mathrm{pr}(\kappa)=e^{a}\star \pi=\pi'.
        \end{equation*} 
        In other words, $\mathrm{pr}_{*}(\tilde{\psi})=\psi$. Hence, the morphism $\mathrm{pr}_{*}$ is a fibration of groupoids.
    \end{proof}
\end{pro}

\begin{lem}\label{twistdg}
If $\alpha\in\mathrm{MC}(\huaL_{V, V'})$, then $\Big(\widetilde{\huaL}_{V, V'}, [\cdot, \cdot]_{\widetilde{\huaL}_{V, V'}}, \dM^{\alpha}=\dM+[\alpha, \cdot]_{\huaL_{V, V'}}\Big)$ is a dgLa, denoted by $\widetilde{\huaL}_{V,V'}^{\alpha}$.
\end{lem}
\begin{proof}
    By Lemma \ref{nonabdgla}, $\big(\widetilde{\huaL}_{V, V'}, [\cdot, \cdot]_{\widetilde{\huaL}_{V, V'}}, \dM\big)$ is an ideal of $\big(\huaL_{V, V'},[\cdot,\cdot]_{\huaL_{V, V'}},d\big)$ as dgLas. Hence, $[\alpha,\xi]_{\huaL_{V, V'}}\in \widetilde{\huaL}_{V, V'}$ for any $\alpha\in\huaL_{V, V'},\xi\in\widetilde{\huaL}_{V, V'}$. Thus, by \Cref{twistdgla}, $\dM^{\alpha}(\widetilde{\huaL}_{V, V'})\subset\widetilde{\huaL}_{V, V'}$, which means that $\Big(\widetilde{\huaL}_{V, V'}, [\cdot, \cdot]_{\widetilde{\huaL}_{V, V'}}, \dM^{\alpha}=\dM+[\alpha, \cdot]_{\huaL_{V, V'}}\Big)$ is a dgLa.
\end{proof}

Using the Deligne groupoid of the twisted dgLa introduced above, we characterize the preimage of a Lie algebra structure on $V'$ under the map $\mathrm{pr}$.
\begin{thm}\label{thm}
Let $\pi\in\mathrm{MC}(\huaL_{V'})$ and $\alpha_0=(0, 0,\pi)\in\mathrm{Obj}(\mathrm{pr}_{*}^{-1}(\pi))\subset\huaL_{V,V'}^{1}$. Then the twisted dgLa $\widetilde{\huaL}_{V,V'}^{\alpha_{0}}$ governs non-abelian extensions of $\g=(V, \mu)$ by $(V', \pi)$. Moreover, we have
$$\mathrm{pr}^{-1}_{*}(\pi)\cong\mathrm{Del}(\widetilde{\huaL}_{V,V'}^{\alpha_{0}}).$$
In particular, when $\pi=0\in\mathrm{MC}(\huaL_{V'})$ is the trivial Lie algebra structure on $V'$, then we have $\mathrm{pr}^{-1}_{*}(0)=\mathrm{Del}(\widetilde{\huaL}_{V,V'})$.
    \begin{proof}
        Let $\pi\in\mathrm{MC}(\huaL_{V'})$ and let $\alpha_0=(0, 0, \pi)\in\MC(\huaL_{V,V'})$ so that $\mathrm{pr}(\alpha_0)=\pi$. In other words, the set of objects 
        $$\mathrm{Obj}\big(\mathrm{pr}_{*}^{-1}(\pi)\big)=\{\kappa\in\MC(\huaL_{V,V'}) : \mathrm{pr}(\kappa)=\pi\}$$ 
        is non-empty. By Lemma \ref{twistdg}, $\widetilde{\huaL}_{V, V'}^{\alpha_{0}}$ is a dgLa. Since $\dM^{\alpha_{0}}(\gamma)=[\mu+\pi, \gamma]$ for all $\gamma\in\widetilde{\huaL}_{V, V'}$, by Proposition \ref{nondgla1}, the twisted dgLa $\widetilde{\huaL}_{V,V'}^{\alpha_{0}}$ governs non-abelian extensions of $\g=(V, \mu)$ by $(V', \pi)$.
        
        For any $\kappa\in\mathrm{Obj}(\mathrm{pr}_{*}^{-1}(\pi))$, we have that $\mathrm{pr}(\kappa-\alpha_{0})=0$, so $\kappa-\alpha_{0}\in\widetilde{\huaL}_{V,V'}$. By the fact that $\kappa$ and $\alpha_{0}$ are Maurer--Cartan elements, we have
        \begin{equation*}
            \dM^{\alpha_{0}}(\kappa-\alpha_0)+\frac{1}{2}[\kappa-\alpha_{0},\kappa-\alpha_{0}]_{\widetilde{\huaL}_{V,V'}}
            =\dM(\kappa)+\frac{1}{2}[\kappa,\kappa]-d(\alpha_{0})
            -\frac{1}{2}[\alpha_{0},\alpha_{0}]=0.
        \end{equation*}
        Hence, we obtain that $\kappa-\alpha_{0}\in\MC(\widetilde{\huaL}_{V,V'}^{\alpha_{0}})$, which means $\mathrm{Obj}(\mathrm{pr}_{*}^{-1}(\pi))-\alpha_{0}\subset\MC(\widetilde{\huaL}_{V,V'}^{\alpha_{0}}).$
        
        Let $u\in\MC(\widetilde{\huaL}_{V,V'}^{\alpha_{0}})$, so $\mathrm{pr}(u+\alpha_{0})=\mathrm{pr}(\alpha_{0})=\pi$. Since
        \begin{eqnarray*}
            d(u+\alpha_{0})+\frac{1}{2}[u+\alpha_{0},u+\alpha_{0}]_{\huaL_{V, V'}}
            =\big(d(u)+\frac{1}{2}[u,u]+[u,\alpha_{0}]\big)+\big(d(\alpha_{0})+\frac{1}{2}[\alpha_{0},\alpha_{0}]\big)=0,
        \end{eqnarray*}
        we obtain $\MC(\widetilde{\huaL}_{V, V'}^{\alpha_{0}})\subset\mathrm{Obj}(\mathrm{pr}_{*}^{-1}(\pi))-\alpha_{0}$. It implies that on the object level, we have 
        $$\mathrm{Obj}\big(\mathrm{pr}_{*}^{-1}(\pi)\big)\cong\mathrm{Obj}\big(\Del(\widetilde{\huaL}_{V, V'}^{\alpha_{0}})\big).$$
       
        On the morphism level, let $\kappa,\kappa'\in\mathrm{Obj}\big( \mathrm{pr}_{*}^{-1}(\pi)\big)$ and $\eta\in\huaL_{V,V'}^{0}=\widetilde{\huaL}_{V, V'}^0$ such that $\kappa'=e^{\eta}\star \kappa$. Denote $u=\kappa-\alpha_{0},\ u'=\kappa'-\alpha_{0}$, then we have
        \begin{equation*}
            u'=u+\sum_{n=0}^{\infty}\frac{\ad^{n}(\eta)}{(n+1)!}\big([\eta,u]-\dM^{\alpha_{0}}(\eta)\big)=e^{\eta}\star u.
        \end{equation*}
        Hence, if $\eta\in\Hom_{\mathrm{pr}_{*}^{-1}(\pi)}(\kappa,\kappa')$, then $\eta\in\Hom_{\Del(\widetilde{\huaL}_{V,V'}^{\alpha_{0}})}(\kappa-\alpha_{0},\kappa'-\alpha_{0})$.
        
        Suppose that $\eta\in\Hom_{\Del(\widetilde{\huaL}_{V,V'}^{\alpha_{0}})}(\kappa-\alpha_{0},\kappa'-\alpha_{0})$, then we obtain 
        \begin{equation*}
\kappa'=\kappa+\sum_{n=0}^{\infty}\frac{\ad^{n}(\eta)}{(n+1)!}\big([\eta,\kappa]-\dM(\eta)\big)=e^{\eta}\star \kappa,
        \end{equation*}
        which implies that  $\eta\in\Hom_{\mathrm{pr}_{*}^{-1}(\pi)}(\kappa,\kappa')$. Thus, $\mathrm{pr}_{*}^{-1}(\pi)\cong\Del(\widetilde{\huaL}_{V,V'}^{\alpha_{0}})$. 

        In particular, when taking $\pi=0\in\mathrm{MC}(\huaL_{V'})$ the trivial Lie algebra structure on $V'$, the isomorphism of the groupoids $\mathrm{pr}^{-1}_{*}(0)$ and $\mathrm{Del}(\widetilde{\huaL}_{V,V'})$ in the above discussion reduces to the identity.
    \end{proof}
\end{thm}
\begin{cor}\label{pro:decom-deligne-gpd}
The Deligne groupoid of $\huaL_{V,V'}$ is the disjoint union of the Deligne groupoids of the twisted dgLas $\widetilde{\huaL}_{V,V'}^{(0,0, \pi)}$ indexed by the Lie algebra structures $\pi$ on $V'$, i.e.
\begin{equation}
    \label{del}\Del\Big(\huaL_{V,V'}\Big)\cong\bigsqcup_{\pi\in\MC(\huaL_{V'})(k)}\Del\Big(\widetilde{\huaL}_{V,V'}^{(0,0,\pi)}\Big).
\end{equation}
\begin{proof}
    Since $\Del(\huaL_{V'})$ is totally disconnected \cite[p. 218]{Br06}, $\mathrm{pr}$ is surjective on $\MC(\huaL_{V,V'})$. By \Cref{thm}, we obtain the desired disjoint union decomposition.
\end{proof}
\end{cor}
\begin{rmk}
    By the fact that $\pi_{0}(\Del(\widetilde{\huaL}_{V,V'}^{(0,0,\pi)}))\cong \mathrm{H}^{2}_{\mathrm{nab}}(V,\mu; V',\pi)$ for all $\pi\in\MC(\huaL_{V'})$ \cite[Proposition 2.2]{Fre14}, taking the set of connected components of the Deligne groupoids in equation \eqref{del}, we recover \Cref{geandnabe}.
\end{rmk}

\section{Obstructions to lifting and extending deformations}\label{sec:obs}

In this section, we study infinitesimal deformations and
obstructions for global extensions. We construct
an obstruction class, the vanishing of which captures the
existence of a lift of an infinitesimal deformation of the
Lie algebra structure on $V'$ to an infinitesimal deformation
of the global extension structure on $V''$.

\subsection{Extending a finite-order deformation}
Let $(V'',p,\g)$ be a global extension datum, where
$\g=(V,\mu)$, and fix a global extension structure $\mu''$.
Let $V'=\ker p$. 

A deformation of $\mu''$ over
$k[t]/(t^{n+1})$ is called an \emph{$n$-th order deformation}.
First-order deformations are also called
\emph{infinitesimal deformations}. We write $D=k[t]/(t^2)$
for the $k$-algebra of dual numbers.

Choose a linear section $s:V\to V''$ of $p$, and let
$(\pi_{V'}^s,\varrho^s,\omega^s)$ be the corresponding
global $2$-cocycle defined in \eqref{2co}. 
Recall the graded spaces $\huaL_{V'}$, $\huaL_{V,V'}$, $\huaL_{V,V'}^{\un}$ and $\widetilde{\huaL}_{V,V'}$ from
\Cref{dglamclie},
\Cref{nondgla2} and \Cref{nondgla1}.
Since $\pi_{V'}^s$ is a Maurer--Cartan element of $\huaL_{V'}$, $\pi^s_{V'}+\varrho^s+\omega^s$ is a Maurer--Cartan element of $\huaL_{V, V'}$, and $\widetilde{\huaL}_{V, V'}$ is an ideal of $\huaL_{V, V'}$. These graded vector spaces give rise to four cochain complexes $(\huaL_{V, V'}, \delta^{V, V'}),\ (\huaL_{V, V'}^{\un}, \delta^{V, V'}),\ (\widetilde{\huaL}_{V, V'}, \widetilde{\delta}^{V, V'})$ and $(\huaL_{V'}, \delta^{V'})$, where 
 $$
 \delta^{V, V'}=[\mu+\pi^s_{V'}+\varrho^s+\omega^s,\cdot], \quad \widetilde{\delta}^{V, V'}=[\mu+\pi^s_{V'}+\varrho^s+\omega^s,\cdot]=\delta^{V, V'}|_{\widetilde{\huaL}_{V, V'}}, \quad \delta^{V'}=[\pi^s_{V'}, \cdot].$$
 Denote by $H^{i}(\huaL_{V, V'})^{\delta},\ H^{i}(\huaL_{V, V'}^{\un})^{\delta},\ H^i(\widetilde{\huaL}_{V, V'})^{\delta}$ and $H^{i}(\huaL_{V'})^{\delta}$ the $i$-th cohomology of the corresponding complexes.
 \begin{rmk}
 Here $\delta^{V,V'}$ obtained by twisting
$d$ introduced in \Cref{nondgla2} by the Maurer--Cartan element
$\pi_{V'}^s+\varrho^s+\omega^s$.
 \end{rmk}
\begin{pro}
    There is a short exact sequence of the cochain complexes:
$$
0\longrightarrow(\widetilde{\huaL}_{V, V‘}, \widetilde{\delta}^{V, V’})\stackrel{\mathrm{\iota}}{\longrightarrow}(\huaL_{V, V‘}, \delta^{V, V’})\stackrel{\mathrm{pr}}{\longrightarrow} (\huaL_{V'}, \delta^{V'})\longrightarrow 0,
$$
where $\iota(\xi)=(\xi, 0)$ and $\mathrm{p}(\xi, f)=f$ for all $f\in \huaL_{V'}$ and $\xi\in \widetilde{\huaL}_{V, V'}$. Consequently,
there is a long exact sequence of the cohomology groups:
\begin{align*}
   0\longrightarrow H^0(\widetilde{\huaL}_{V, V'})^{\delta}&{\longrightarrow} H^0(\huaL_{V, V'})^{\delta}{\longrightarrow} 0\longrightarrow  H^{1}(\widetilde{\huaL}_{V, V'})^{\delta} \\
&\longrightarrow H^1(\huaL_{V, V'})^{\delta}{\longrightarrow} H^1(\huaL_{V'})^{\delta}\longrightarrow H^2(\widetilde{\huaL}_{V, V'})^{\delta}\longrightarrow\cdots. 
\end{align*}

\end{pro}
 \begin{proof}
The short exact sequence follows from a straightforward computation.
 \end{proof}
In fact, various deformations information about global extensions can be read from the cohomology of the dgLa, as in the classical setting. A typical example is that the first cohomology can be used to characterize infinitesimal of $\mu''$.
\begin{pro}\label{infidef}
    Equivalence (resp. unframed equivalence) classes of infinitesimal deformations of the global extension structure $\mu''$ are classified by the first cohomology of $(\huaL_{V, V'}, \delta^{V, V'})$ (resp. $(\huaL_{V, V'}^{\un}, \delta^{V, V'})$). In other words, there are bijections 
    \begin{equation*}
               \mathcal{D}_{\mu''}\big(k[t]/(t^{2})\big)\cong H^{1}(\huaL_{V,V'})^{\delta},\quad \mathcal{D}_{\mu''}^{\un}\big(k[t]/(t^{2})\big)\cong H^{1}(\huaL_{V,V'}^{\un})^{\delta},
    \end{equation*}
    where $\mathcal{D}_{\mu''}=\pi_{0}\mathfrak{Def}_{\mathcal{E}_p, [\mu^{''}]}$ and $\mathcal{D}_{\mu''}^{\un}=\pi_{0}\mathfrak{Def}_{\mathcal{E}^{\un}_{p}, [\mu^{''}]}$ are set-valued deformation functors induced by $\mathfrak{Def}_{\mathcal{E}_p, [\mu^{''}]}$ and $\mathfrak{Def}_{\mathcal{E}^{\un}_{p}, [\mu^{''}]}$ in \Cref{4} respectively.
    \begin{proof}
     This follows from \Cref{dgladefglobal}, \Cref{unframed-formal-moduli} and \cite[Lemma 6.4.1]{Man22}.
    \end{proof}
\end{pro}
For a fixed algebraic structure, the question of whether and under what conditions an $n$-th order deformation can be extended to an $(n+1)$-th order deformation has also been extensively studied in the literature. These problems can be transformed into a cohomological problem, which is now known as the obstruction theory.

\begin{pro}\label{pro:n-order-obs}
  Let $\mu''$ be a fixed global extension structure on $V''$, and $\mu_{n}{''}$ an $n$-th order deformation of $\mu''$. Then $\mu_{n}{''}$ is extensible to an $(n+1)$-th order deformation if and only if the obstruction class 
    \begin{equation*}
       \mathrm{Ob}_{n}(\mu''_{n})\coloneqq\frac{1}{2}\Big[\sum_{i+j=n+1}[\eta_i,\eta_j]\Big]\in H^{2}(\huaL_{V,V'})^{\delta}
    \end{equation*} vanishes. Here, $\mu''=\mu+\sum_{i=1}^{n}t^i\eta_i,$ and $ \eta_i=\pi_{V',i}+\varrho_i+\omega_i$ are the coefficients of $t^i$ in the expansion of the corresponding global 2-cocycle of $\mu''_{n}-\mu''$ as a polynomial in $t$.
    \begin{proof}
        This follows from \Cref{dgladefglobal} and \cite[\S 3]{Man99}.
    \end{proof}
\end{pro}  

\subsection{Lifting a deformation of the Lie bracket on the kernel}
Fix a linear section $s$ of $p$. For convenience, we will omit the superscript $s$ if there is no ambiguity. By \Cref{defglobal2co},
an infinitesimal deformation $\mu''_D$ of $\mu''$ corresponds
to a global $2$-cocycle
$$
 (\pi_{V'}^D,\varrho^D,\omega^D) =
 (\pi_{V'}+t\pi_{V'}',\varrho+t\varrho',\omega+t\omega'),
$$
where
$$
 \pi_{V'}'\in\Hom(\wedge^2V',V'),\quad
 \varrho'\in\Hom(V\otimes V',V'),\quad
 \omega'\in\Hom(\wedge^2V,V').
$$
In particular, $\pi_{V'}^D$ is an infinitesimal deformation
of the Lie bracket $\pi_{V'}$ on the kernel.

Conversely, given such a deformation
$\pi_{V'}^D=\pi_{V'}+t\pi_{V'}'$, we ask whether it
extends to an infinitesimal deformation of the global
cocycle $(\pi_{V'},\varrho,\omega)$, and determine
the obstruction to finding suitable $\varrho'$ and $\omega'$ such that 
\begin{equation*}
   (\pi_{V'}^{D}, \varrho^{D}, \omega^{D})\coloneqq (\pi_{V'}+t\pi_{V'}',\ \varrho+t\varrho',\ \omega+t\omega')
\end{equation*}
forms an infinitesimal deformation of the given global 2-cocycle $(\pi_{V'},\ \varrho,\ \omega)$.

Recall that if $\pi_{V'} + t\pi_{V'}'$ is an infinitesimal deformation of $\pi_{V'}$, according to \cite[\S 2]{NR64}, $\pi_{V'}'$ is a cocycle in the corresponding Chevalley-Eilenberg complex. Thus, $\delta^{V'}(\pi_{V'}')=0$.

\begin{pro}\label{liftob}
There exists a pair $(\varrho^{'},\omega^{'})\in\Hom(V\otimes V', V')\times\Hom(\wedge^2V, V')$ such that
\begin{equation*}
    (\pi_{V'}^{D}, \varrho^{D}, \omega^{D})\coloneqq (\pi_{V'}+t\pi_{V'}',\ \varrho+t\varrho',\ \omega+t\omega')
\end{equation*}
is an infinitesimal deformation of the given global 2-cocycle $(\pi_{V'},\ \varrho,\ \omega)$ if and only if the cohomology class $\mathrm{Ob}_{\pi_{V'}'}(\pi_{V'}, \varrho, \omega)\in H^{2}(\widetilde{\huaL}_{V, V'})^{\delta}$, defined by
$$
\mathrm{Ob}_{\pi_{V'}'}(\pi_{V'}, \varrho, \omega)=\Big[[\varrho+\omega, \pi_{V'}']\Big],
$$
vanishes. 
\end{pro}
\begin{proof}
 By \eqref{defibra}, we know that $[\mu, \pi_{V'}']=0.$ Since $\delta^{V'}(\pi_{V'}')=[\pi_{V'}, \pi_{V'}']=0$, we have $[\varrho+\omega, \pi_{V'}']=\widetilde{\delta}^{V, V'}(\pi_{V'}')$. Thus, there is a cohomology class $\mathrm{Ob}_{\pi_{V'}'}(\pi_{V'}, \varrho, \omega)\in H^{2}(\widetilde{\huaL}_{V, V'})^{\delta}$.

 If $\mathrm{Ob}_{\pi_{V'}'}(\pi_{V'}, \varrho, \omega)=0$, then there exist $\varrho' \in \text{Hom}(V\otimes V', V')$ and $\omega' \in \text{Hom}(\wedge^2V, V')$ such that 
\begin{equation*}
    [\varrho+\omega, \pi_{V'}']=\widetilde{\delta}^{V, V'}(-\varrho'-\omega'),
\end{equation*}
which is equivalent to the following equations 
\begin{equation*}
    [\varrho, \pi_{V'}']=-[\pi_{V'}, \varrho'],\quad [\omega, \pi_{V'}']=-[\mu, \varrho']-[\varrho, \varrho']-[\pi_{V'}, \omega'],\quad [\mu, \omega']=-[\varrho, \omega']-[\varrho', \omega].
\end{equation*}
Therefore, $(\omega+t\omega', \varrho+t\varrho',  \pi_{V'} + t\pi_{V'}')$ satisfies the equations \eqref{gc4}, \eqref{gc5}, and \eqref{gc6}. Define $\varrho^{D}=\varrho+t\varrho'$ and $\omega^{D}=\omega+t\omega'$. Since
\begin{equation*}
    \varrho^{D}(x\otimes 1)(u\otimes 1)=\varrho(x)u\otimes 1+\varrho'(x)u\otimes t, \quad \omega^{D}(x\otimes 1, y\otimes 1)=\omega(x, y)\otimes 1+\omega'(x, y)\otimes t,
\end{equation*}
it follows that $(\pi_{V'}^{D}, \varrho^{D}, \omega^{D})$ is a global $2$-cocycle, which induces an infinitesimal deformation $\mu''_{D}$ of $\mu''$.

If $(\pi_{V'}^{D}, \varrho^{D}, \omega^{D})$ is a global $2$-cocycle, then there exists $$(\varrho^{'},\omega^{'})\in\Hom(V\otimes V', V')\times\Hom(\wedge^2V, V')$$
such that
\begin{equation*}
    \varrho^{D}=\varrho+t\varrho', \quad \omega^{D}=\omega+t\omega'.
\end{equation*}
The global $2$-cocycle conditions imply the following identities:
\begin{align*}
~~[\varrho, \pi_{V'}']&=-[\pi_{V'}, \varrho'],\\
~~[\omega, \pi_{V'}']&=-[\mu, \varrho']-[\varrho, \varrho']-[\pi_{V'}, \omega'],\\
~~[\mu, \omega']&=-[\varrho, \omega']-[\varrho', \omega].
\end{align*}
Consequently, these relations can be unified as
\begin{equation*}
    [\varrho+\omega, \pi_{V'}']=\widetilde{\delta}^{V, V'}(-\varrho'-\omega'),
\end{equation*}
which implies that $\mathrm{Ob}_{\pi_{V'}'}(\pi_{V'}, \varrho, \omega)=0$.  
\end{proof}
\begin{ex}
We explicate the dgLa description for \Cref{ex11}. We assume that $\g=(V, \mu)$ has even-dimension $2n$.
Let $\{y_1, \cdots,y_n, z_1,\cdots, z_{n}\}$ be a basis of $V$. Define $\varrho$ and $\omega$ by
$$\varrho=0, \quad \omega(y_i, z_j)=\delta_{ij}, \quad \omega(y_k, y_l)=\omega(z_p, z_q)=0.$$
    Then the global extension $(k\oplus V, [\cdot, \cdot]_{\omega})$ is the $(2n+1)$-dimensional Heisenberg Lie algebra, where 
    $$
    [( r_1,x_1), (r_2,x_2)]_{\omega}=(\omega(x_1, x_2),0), \quad \forall (r_i, x_i)\in k\oplus V,\quad i=1,2.
    $$
     By Theorem \ref{dgladefglobal}, the dgLa 
    $$\huaL_{V,k}^{\omega}=\Big(\huaL_{V,k},[\cdot,\cdot]_{\huaL_{V,k}},d=[\omega,\cdot]_{\huaL_{V, k}}\Big)$$
    controls the deformation of $(k\oplus V, [\cdot, \cdot]_{\omega}, V, p)$, 
    where
    \begin{equation*}
        \huaL_{V, k}^{i}=
        \begin{cases}
        0, & i\leq-1;\\
        \Hom(V, k), & i=0;\\
        \Hom(\wedge^{i+1}V,k)\oplus\Hom(\wedge^{i}V\otimes k,k),\ & i\geq 1.
        \end{cases}
     \end{equation*}
    By Proposition \ref{infidef}, the set of equivalence classes of infinitesimal deformations of $(k\oplus V, [\cdot, \cdot]_{\omega}, V, p)$ is 
    \begin{equation*}
        \huaD_{[\cdot, \cdot]_{\omega}}\big(k[t]/(t^{2})\big)\cong \mathrm{H}^{1}(\huaL_{V, k}^{\omega}),
    \end{equation*} 
    where $\mathcal{D}_{[\cdot, \cdot]_{\omega}}$ is the set-valued deformation functor induced by $\mathfrak{Def}_{\huaE_p, [[\cdot, \cdot]_{\omega}]}$.
    By \eqref{defibra}, we have $$df=[\omega, f]_{\huaL_{V, k}}=0, \quad\forall f\in\Hom(V, k), \quad\text{and} \quad[\omega, \Hom(\wedge^{2}V, \ k)]_{\huaL_{V, k}}=0.$$ For $\varpi\in\Hom(V\otimes k, k)\cong V^{*}$, since the differential is $0$ on $\Hom(\wedge^{2}V,k)$, a direct computation shows that
 
   \begin{equation*}
       H^{1}(\huaL_{V, k}^{\omega})\cong
\wedge^{2}V^{*}\oplus\Ker(V^{*}\overset{\varpi\mapsto\varpi\wedge\omega}{\longrightarrow}\wedge^{3}V^{*})
    \end{equation*} 
Notice that $\wedge^{3}V^{*}=0$ for $n=1$, and the map $\varpi\mapsto\varpi\wedge\omega$ is injective for $n=2$, which implies that 
 \begin{align*}
H^{1}(\huaL_{V, k}^{\omega})\cong\begin{cases}
           \wedge^{2}V^{*}\oplus V^{*}, & n=1,\\
           \wedge^{2}V^{*}, & n\geq 2.
         \end{cases}
    \end{align*} 
 
\end{ex}
The following example shows that the obstruction in \Cref{liftob} can be nonzero.
\begin{ex}
Let $V=\rmspan_{k}\{x\}$ be a one-dimensional abelian Lie algebra with $\mu=0$, and let $V'=\rmspan_{k}\{u,v\}$. Denote by $V''=V'\oplus V\overset{p}{\lon}V$ the projection to the second summand, $s:V\lon V''$ the natural embedding. The induced global $2$-cocycle is 
\begin{equation*}
    \pi_{V'}=0,\quad
    \varrho(x)=\operatorname{id}_{V'},\quad
    \omega=0.
\end{equation*}
It determines a global extension on $V''=V'\oplus V$, whose
Lie bracket is given by
\begin{equation*}
    \mu''(x, u)=u,\quad \mu''(x, v)=v.
\end{equation*}
Let $D=k[t]/(t^2)$. Define
\begin{equation*}
    \pi'_{V'}\in\operatorname{Hom}(\wedge^2 V',V'),
    \quad \pi'_{V'}(u,v)=v.
\end{equation*}
Since $\pi_{V'}=0$, we have $[\pi_{V'},\pi'_{V'}]=0.$ Thus
$$
    \pi^D_{V'}=\pi_{V'}+t\pi'_{V'}=t\pi'_{V'}
$$
is an infinitesimal deformation of the Lie bracket on the kernel. We compute the obstruction in \Cref{liftob}. Set
\begin{equation*}
    O=[\varrho+\omega,\pi'_{V'}]
      =[\varrho,\pi'_{V'}].
\end{equation*}
By direct calculation, we obtain
\begin{align*}
    O(x,u,v)=\pi'_{V'}(\varrho(x)u,v)+\pi'_{V'}(u,\varrho(x)v) -\varrho(x)\pi'_{V'}(u,v)=v+v-v=v.
\end{align*}
In particular, $O\neq 0$.

Moreover, since $\dim V=1$, we have
\begin{equation*}
     \widetilde L^1_{V,V'} =\operatorname{Hom}(V\otimes V',V'),\quad
    \widetilde L^2_{V,V'}=\operatorname{Hom}(V\otimes\wedge^2 V',V'),\quad
    \widetilde L^3_{V,V'}=0.
\end{equation*}
The differential on $\widetilde L_{V,V'}$ is
$$
    \widetilde\delta^{V,V'}
      =[\mu+\pi_{V'}+\varrho+\omega,-]
      =[\varrho,-].
$$
For every
$\varrho'\in\operatorname{Hom}(V\otimes V',V')$, the bracket
$[\varrho,\varrho']$ belongs to
$\operatorname{Hom}(\wedge^2 V\otimes V',V')=0$.
Consequently,
$$
    \widetilde\delta^{V,V'}
       \bigl(\widetilde L^1_{V,V'}\bigr)=0.
$$
Hence $O$ is a cocycle which is not a coboundary, and
\begin{equation*}
     \operatorname{Ob}_{\pi'_{V'}}(\pi_{V'},\varrho,\omega)=[O]\neq 0\quad\text{in }H^2(\widetilde L_{V,V'}).
\end{equation*}
By \Cref{liftob}, the deformation $\pi^D_{V'}$ of the kernel
does not lift to an infinitesimal deformation of the given
global extension.

This can also be seen directly. Any such lift would satisfy
$$    \varrho^D(x)=\operatorname{id}_{V'}+t\varrho'(x),
    \quad \omega^D=0.
$$
The requirement that $\varrho^D(x)$ is a derivation of
$\pi^D_{V'}$, evaluated on $u,v$, would give
$$
    \varrho^D(x)\pi^D_{V'}(u,v)=
    \pi^D_{V'}(\varrho^D(x)u,v)
    +\pi^D_{V'}(u,\varrho^D(x)v),
$$
that is, $tv=2tv$, a contradiction in $V'\otimes D$.
\end{ex}

\section{An application: contraction from $\mathfrak{g}\rtimes_{\mathrm{ad}} \mathfrak{g}$ to $V_{\mathrm{ab}}\rtimes_{\mathrm{ad}}\g$}\label{sec:contraction}

\begin{defi}\cite[\S1.1]{FF99}\label{def:algebraic-deformation}
Let $A$ be an augmented commutative $\K$-algebra. An \emph{algebraic deformation} of a Lie algebra $\mathfrak{g}=(V,\nu)$ over $A$ is an $A$-bilinear Lie algebra structure on the free module $V\otimes A$ whose reduction along the augmentation is $\nu$. 
\end{defi}

\begin{rmk}
This definition does not require $A$ to be local or Artinian. A closely related notion is the algebraic contraction. An \emph{algebraic contraction} of $\nu_1$ to $\nu_0$ is a family over $\K[t]$, augmented by $t\mapsto0$, whose special fiber is $\nu_0$ and whose restriction to $\mathbb G_m=\mathbb A^1\setminus\{0\}$ is isomorphic to the constant family with fiber $\nu_1$.
\end{rmk}
We will also deal with the case when $(V'',\mu'')$ is a global extension structure. In this case, the notions of deformation and contraction are defined similarly, keeping track of the extension structure. 
\begin{defi}\label{framecontraction}
We call the contraction \emph{kernel-framed} if the isomorphism over $\mathbb G_m$ fixes $V'$ pointwise, and \emph{kernel-unframed} if $V'$ is only stable under the isomorphism over $\mathbb G_m$. 
\end{defi}

Let $\g=(V,\mu)$, and let $V'$ be a second copy of its underlying vector space. Denote the transported bracket on $V'$ by $\mu_{V'}(u,v)=[u,v]$ and the adjoint action by $\ad_x(u)=[x,u]$. Set $A=\K[t]$, augmented by $t\mapsto0$, and write $V'_{A}=V'\otimes A$ and $\g_A=(V\otimes A,\mu\otimes1_A)$.

\begin{pro}\label{prop:adjoint-family}
The formula
\begin{equation}\label{eq:adjoint-family}
[(u,x),(v,y)]_t
=\bigl([x,v]-[y,u]+t[u,v],[x,y]\bigr)
\end{equation}
defines an $A$-Lie algebra structure on $(V'\oplus V)\otimes A$. The projection to $\g_A$ is a Lie algebra homomorphism. With respect to the augmentation $t\mapsto0$, it is an algebraic deformation of $V'_{\mathrm{ab}}\rtimes_{\ad}\g$.
\end{pro}

\begin{proof}
Endow $V'_{A}$ with the bracket $[u,v]_t=t[u,v]$. Its Jacobiator is $t^2$ times the Jacobiator of $\g$, so it vanishes. The adjoint action satisfies
\begin{equation*}
    [x,[u,v]_t]=[[x,u],v]_t+[u,[x,v]]_t
\end{equation*}
and $[\ad_x,\ad_y]=\ad_{[x,y]}$. Therefore, we get an action of $\g_A$ on $(V'_{A},[\cdot,\cdot]_t)$ by derivations. Moreover, the corresponding Lie algebra structure is given by \eqref{eq:adjoint-family}, which specializes to $V'_{\mathrm{ab}}\rtimes_{\ad}\g$ when $t=0$.
\end{proof}

\begin{pro}\label{prop:adjoint-fibers}
The family in \Cref{prop:adjoint-family} is flat. Its fiber at $t=0$ is $V'_{\mathrm{ab}}\rtimes_{\ad}\g$, where $V'_{\mathrm{ab}}$ is abelian. Its fiber at $t=1$ is $\g\rtimes_{\ad}\g$, where the second factor carries its original bracket. For every $a\in\K^\times$, the map
\begin{equation*}
    \theta_a:(V'',[\cdot,\cdot]_a)\longrightarrow(V'',[\cdot,\cdot]_1),
\quad \theta_a(u,x)=(au,x),
\end{equation*}
is a kernel-unframed global extension isomorphism. Consequently,
$\g\rtimes_{\ad}\g$ contracts to
$V'_{\mathrm{ab}}\rtimes_{\ad}\g$ at $t=0$ in the kernel-unframed sense of \Cref{framecontraction}.

The punctured line $\mathbb A^1\setminus\{0\}=\mathbb G_m$ is a $\mathbb G_m$-torsor under multiplication. This action lifts to the punctured family and extends to the $\mathbb A^1$-family by
\begin{equation}\label{eq:gm-action}
\lambda\cdot(t;u,x)=(\lambda t;\lambda^{-1}u,x),\quad \forall \lambda\in \mathbb G_m,t\in \mathbb A^1, u\in V', x\in V.
\end{equation}
\end{pro}

\begin{proof}
Since the underlying $A$-module is free, the family is flat. By substituting $t=0$ and $t=1$ in \eqref{eq:adjoint-family}, we get the descriptions of the corresponding fibers. If $a\in\K^\times$, then
\begin{align*}
\theta_a([(u,x),(v,y)]_a)
&=\bigl(a[x,v]-a[y,u]+a^2[u,v],[x,y]\bigr)\\
&=[\theta_a(u,x),\theta_a(v,y)]_1.
\end{align*}
In particular, the restriction $a\Id_{V'}$ is a Lie algebra isomorphism
from $(V',a\mu_{V'})$ to $(V',\mu_{V'})$, since
\begin{equation*}
    a\Id_{V'}\bigl(a[u,v]\bigr)=a^2[u,v]=[au,av].
\end{equation*}
We remark that the map $a\Id_{V'}$ does not belong to $\operatorname{Aut}(V',\mu_{V'})$, and it is a Lie algebra isomorphism only after we rescale the Lie bracket on $V'$. Hence, the restriction of the family to $\mathbb G_m$ is trivialized by
\begin{equation*}
    (t;u,x)\longmapsto(t;tu,x).
\end{equation*}
Moreover, the kernel at $t=0$ is abelian. This produces a particular example of the usual In\"on\"u--Wigner contraction construction as in \cite{PY}.

It remains to check that the action actually lifts. Let
$\Theta_\lambda(u,x)=(\lambda^{-1}u,x)$. Then
\begin{align*}
\Theta_\lambda([(u,x),(v,y)]_t)
&=\bigl(\lambda^{-1}[x,v]-\lambda^{-1}[y,u]
  +\lambda^{-1}t[u,v],[x,y]\bigr)\\
&=[\Theta_\lambda(u,x),\Theta_\lambda(v,y)]_{\lambda t}.
\end{align*}
Hence, the action in \eqref{eq:gm-action} carries the fiber at $t$ isomorphically to the fiber at $\lambda t$. Moreover, it is direct to check
\begin{equation*}
    \Theta_{\lambda\eta}=\Theta_\lambda\circ\Theta_\eta,
\quad
(\lambda\eta)t=\lambda(\eta t),
\end{equation*}
so that these fiberwise Lie algebra isomorphisms define a $\mathbb G_m$-action on the punctured family.
\end{proof}

\begin{pro}\label{prop:adjoint-mc}
Let $s_{\mathrm{ab}}=\ad\in\MC(L_{V,V'})$ be the Maurer--Cartan element for the special fiber at $t=0$. Relative to the twist $L_{V,V'}^{s_{\mathrm{ab}}}$, the family is
\begin{equation*}
    \beta(t)=t\mu_{V'}.
\end{equation*}
For every $n\geq1$,
\begin{equation*}
    t\mu_{V'}\in\MC\!\left(L_{V,V'}^{s_{\mathrm{ab}}}\otimes
(t)\K[t]/(t^{n+1})\right).
\end{equation*}
If $\g$ is non-abelian, the first-order classes
\begin{equation*}
    [\mu_{V'}]\in H^1(L_{V,V'}^{s_{\mathrm{ab}}}),
\quad
[\mu_{V'}]\in H^1(L_{V,V'}^{\un,s_{\mathrm{ab}}})
\end{equation*}
are nonzero. Thus, the formal germ at the degenerate fiber is nontrivial.
\end{pro}

\begin{proof}
It is straightforward to see that the absolute Maurer--Cartan element of \eqref{eq:adjoint-family} is
\begin{equation*}
    s(t)=\ad+t\mu_{V'}=s_{\mathrm{ab}}+t\mu_{V'}.
\end{equation*}
 Since the adjoint action is by derivations, together with the Jacobi identity, we have 
 $$[\ad,\mu_{V'}]_{\NR}=0,\quad [\mu_{V'},\mu_{V'}]_{\NR}=0.$$

The projection $\mathrm{pr}:L_{V,V'}\to L_{V'}$ in \Cref{nonabdgla} sends $d^{s_{\mathrm{ab}}}\tau$ to zero for every $\tau\in L_{V,V'}^0$. It sends the cocycle $\mu_{V'}$ to the nonzero bracket $\mu_{V'}$. Hence $\mu_{V'}$ cannot be a boundary when $\g$ is non-abelian.

The same argument applies to the enlarged degree-zero term. If
$\xi\in L_{V,V'}^{\un, 0}=\Hom(V'',V')$, then the restriction of
$d^{s_{\mathrm{ab}}}\xi=[\mu+s_{\mathrm{ab}},\xi]_{\NR}$ to
$V'\wedge V'$ is zero, because the Lie bracket on the kernel of the special fiber is zero. The restriction of $\mu_{V'}$ to $V'\wedge V'$ is the nonzero original bracket. Therefore, $\mu_{V'}$ is not a boundary in $L_{V,V'}^{\un ,s_{\mathrm{ab}}}$ either.

Notice that the isomorphism $\theta_a$ changes the framing on the kernel, so that it is an arrow only in the kernel-unframed groupoid. After inverting $t$, the formal family becomes isomorphic to the constant family with fiber $\g\rtimes_{\ad}\g$ via the same change.
\end{proof}

The following result is a direct consequence of the kernel-unframed global extension formulism in \Cref{thm:unframed-stack}.

\begin{pro}\label{prop:adjoint-global-groupoid}
Assume that $\g$ is non-abelian. Let
\begin{equation*}
    c:\mathbb A^1\longrightarrow \MC(L_{V,V'}), \quad
c(t)=s(t)=s_{\mathrm{ab}}+t\mu_{V'}.
\end{equation*}
Then the following statements hold.
\begin{enumerate}[(i)]

\item For $\lambda\in\mathbb G_m$, let
\begin{equation*}
    \Theta_\lambda=\lambda^{-1}\Id_{V'}\oplus\Id_V\in H_{V,V'}.
\end{equation*}
This is a kernel-unframed global extension isomorphism from the fiber at
$t$ to the fiber at $\lambda t$, and transport of brackets satisfies
\begin{equation*}
    \Theta_\lambda\cdot c(t)=c(\lambda t).
\end{equation*}
Consequently, the restriction of the family to $\mathbb G_m$ is constant in the kernel-unframed stack $\mathfrak{MC}(L_{V, V'}^{\un})$, while its value at $0$ is
the contraction limit. Thus, $c$ is an algebraic deformation of global extensions and a kernel-unframed algebraic contraction.

\item Put $\nu_{\mathrm{ab}}=\mu+s_{\mathrm{ab}}$. In the kernel-framed
and kernel-unframed based formal deformation problems, the germ of $c$ is the Maurer--Cartan family $t\mu_{V'}$ in, respectively,
$L_{V,V'}^{s_{\mathrm{ab}}}$ and $L_{V,V'}^{\un,s_{\mathrm{ab}}}$. Its tangent vectors are the nonzero classes
\begin{equation*}
    [\mu_{V'}]\in H^1(L_{V,V'}^{s_{\mathrm{ab}}}), \quad
[\mu_{V'}]\in H^1(L_{V,V'}^{\un,s_{\mathrm{ab}}}).
\end{equation*}
Moreover, the rescaling trivializes the punctured algebraic family.
\end{enumerate}
\end{pro}

\begin{proof}
On the kernel, the rescaling has the identity
\begin{equation*}
    \lambda^{-1}\Id_{V'}\bigl(t[u,v]\bigr)=\lambda^{-1}t[u,v]
=[\lambda^{-1}u,\lambda^{-1}v]_{\lambda t}.
\end{equation*}
Together with the compatibility of the adjoint action, this gives the
kernel-unframed isomorphism
$\Theta_\lambda:(V'',[\cdot,\cdot]_t)\to(V'',[\cdot,\cdot]_{\lambda t})$ and therefore
$\Theta_\lambda\cdot c(t)=c(\lambda t)$. Notice that
$\lambda^{-1}\Id_{V'}$ is an isomorphism from $(V',t\mu_{V'})$ to
$(V',\lambda t\mu_{V'})$. It follows that every point $c(a)$ with $a\ne0$ lies in the torus orbit of $c(1)$. The trivialization over $\mathbb G_m$ is given by $\theta_t= t\Id_{V'}\oplus\Id_V$.

The last statement follows from \Cref{prop:adjoint-mc} together with the two equivalences
\begin{equation*}
    \mathfrak{Def}_{\mathcal E_p,[\nu_{\mathrm{ab}}]}(A)\simeq
\Del(L_{V,V'}^{s_{\mathrm{ab}}}\otimes\m_A),
\quad
\mathfrak{Def}_{\mathcal E^{\un}_p,[\nu_{\mathrm{ab}}]}(A)\simeq
\Del(L_{V,V'}^{\un,s_{\mathrm{ab}}}\otimes\m_A)
\end{equation*}
in \Cref{dgla2co}, \Cref{defglobal2co} and
\Cref{unframed-formal-moduli}.
\end{proof}

\bibliographystyle{amsplain}
\bibliography{referenceglobal}

 \end{document}